\documentclass[twoside, 11pt]{article}

\usepackage{amssymb}
\usepackage{amsfonts}
\usepackage{amsthm}
\usepackage{amsmath}
\usepackage{amscd}

\begin{document}

\pagestyle{myheadings}
\markboth{\c S.A. Basarab\/}
{\c S.A. Basarab\/}
\title{{\huge Embedding theorems for 
actions on generalized trees, I\/}}

\author{\c Serban A. Basarab
\thanks{The author gratefully thanks Professor 
Thomas M\" uller for sending
him his nice joint work with Ian Chiswell 
{\em ``Embedding theorems for tree-free groups''},
which motivated and inspired the present work.}\\
Institute of Mathematics "Simion Stoilow" of the
Romanian Academy \\
P.O. Box 1--764\\
RO -- 014700 Bucharest, ROMANIA\\
\texttt{{\itshape e-mail}: Serban.Basarab@imar.ro}\\ 
\\
Dedicated to the memory of my father,\\
the painter Alexandru Bassarab (1907-1941)}
\date{} 

\maketitle
\vspace{8mm}

\newtheorem{te}{Theorem}[section]
\newtheorem{ex}[te]{Example}
\newtheorem{pr}[te]{Proposition}
\newtheorem{rem}[te]{Remark}
\newtheorem{rems}[te]{Remarks}
\newtheorem{co}[te]{Corollary}
\newtheorem{lem}[te]{Lemma}
\newtheorem{prob}[te]{Problem}
\newtheorem{exs}[te]{Examples}
\newtheorem{defs}[te]{Definitions}
\newtheorem{de}[te]{Definition}
\newtheorem{que}[te]{Question}
\newtheorem{aus}[te]{}


\newenvironment{aussage}%
{\renewcommand{\theequation}{\alph{equation}}\setcounter{equation}{0}
\begin{aussage*}}{\end{aussage*}}

\def\R{\mathbb{R}}
\def\bbbm{{\rm I\!M}}
\def\N{\mathbb{N}}
\def\F{\mathbb{F}}
\def\H{\mathbb{H}}
\def\I{\mathbb{I}}
\def\K{\mathbb{K}}
\def\P{\mathbb{P}}
\def\D{\mathbb{D}}
\def\Z{\mathbb{Z}}
\def\C{\mathbb{C}}
\def\T{\mathbb{T}}
\def\L{\mathbb{L}}
\def\X{\mathbb{X}}
\def\G{\mathbb{G}}

\def\cA{\mathcal{A}}
\def\cB{\mathcal{B}}
\def\cC{\mathcal{C}}
\def\cD{\mathcal{D}}
\def\cE{\mathcal{E}}
\def\cF{\mathcal{F}}
\def\cG{\mathcal{G}}
\def\cH{\mathcal{H}}
\def\cI{\mathcal{I}}
\def\cJ{\mathcal{J}}
\def\cK{\mathcal{K}}
\def\cL{\mathcal{L}}
\def\cM{\mathcal{M}}
\def\cN{\mathcal{N}}
\def\cO{\mathcal{O}}
\def\cP{\mathcal{P}}
\def\cQ{\mathcal{Q}}
\def\cR{\mathcal{R}}
\def\cS{\mathcal{S}}
\def\cT{\mathcal{T}}
\def\cU{\mathcal{U}}
\def\cV{\mathcal{V}}
\def\cW{\mathcal{W}}
\def\cX{\mathcal{X}}
\def\cY{\mathcal{Y}}
\def\cZ{\mathcal{Z}}


\newcommand{\notdiv}{{\not{|}\,}}
\newcommand{\ctg}{{\rm ctg\,}}
\newcommand{\sh}{{\rm sh\,}}
\newcommand{\ch}{{\rm ch\,}}
\newcommand{\Spin}{{\rm Spin\,}}
\newcommand{\spin}{{\rm spin\,}}
\newcommand{\Spek}{{\rm Spek\,}}
\newcommand{\Ker}{{\rm Ker\,}}
\newcommand{\cont}{\mbox{\rm cont\,}}
\newcommand{\rank}{\mbox{\rm rank\,}}
\newcommand{\codim}{\mbox{\rm codim\,}}
\newcommand{\ssqrt}[2]{\sqrt[\scriptstyle{#1}]{#2}}
\newcommand{\dbigcup}{\displaystyle\bigcup}
\newcommand{\dprod}{\displaystyle\prod}
\newcommand{\dbigcap}{\displaystyle\bigcap}
\newcommand{\dinf}{\displaystyle\inf}
\newcommand{\dint}{\displaystyle\int}
\newcommand{\dsum}{\displaystyle\sum}
\newcommand{\dbigoplus}{\displaystyle{\bigoplus_{}}}
\newcommand{\comb}[2]{\left(\begin{array}{c}#1\\#2\end{array}\right)}
\def\be{\begin{equation}}
\def\ee{\end{equation}}
\newcommand{\rs}{respectiv }
\newcommand{\dleq}{\displaystyle\leq}
\newcommand{\dgeq}{\displaystyle\geq}

\def\bp{\begin{proof}}
\def\ep{\end{proof}}
\def\hf{\hfill $\square$}
\def\ben{\begin{enumerate}}
\def\een{\end{enumerate}}
\def\ba{\begin{eqnarray*}}
\def\ea{\end{eqnarray*}}
\def\del{\delta}
\def\ve{\varepsilon}
\def\ls{\leqslant}
\def\gs{\geqslant}
\def\lra{\longrightarrow}
\def\Llra{\Longleftrightarrow}
\def\Lra{\Longrightarrow}
\def\p{\perp}
\def\wp{\,\widehat{1/p}}
\def\w4{\,\widehat{1/4}}
\newcommand{\wh}{\widehat}
\newcommand{\la}{\langle}
\newcommand{\ra}{\rangle}
\newcommand{\wt}{\widetilde}
\newcommand{\sm}{\setminus}
\newcommand{\sse}{\subseteq}
\newcommand{\es}{\varnothing} 
\newcommand{\fG}{\mathfrak{G}}
\newcommand{\fB}{\mathfrak{B}}
\newcommand{\fZ}{\mathfrak{Z}}
\newcommand{\fK}{\mathfrak{K}}
\newcommand{\fQ}{\mathfrak{Q}}
\newcommand{\fL}{\mathfrak{L}}
\newcommand{\fC}{\mathfrak{C}}
\newcommand{\fS}{\mathfrak{S}}
\newcommand{\fT}{\mathfrak{T}}
\newcommand{\fM}{\mathfrak{M}}
\newcommand{\fU}{\mathfrak{U}}
\newcommand{\fP}{\mathfrak{P}}
\newcommand{\fH}{\mathfrak{H}}
\newcommand{\fX}{\mathfrak{X}}
\newcommand{\fn}{\mathfrak{n}}
\newcommand{\fm}{\mathfrak{m}}
\newcommand{\fp}{\mathfrak{p}}
\newcommand{\fq}{\mathfrak{q}}
\newcommand{\f}{\frac}
\newcommand{\q}{\quad}
\newcommand{\n}{\vartriangleleft}

\newcommand{\al}{h }
\newcommand{\De}{\Delta}
\newcommand{\eps}{\varepsilon}
\newcommand{\gam}{\gamma}
\newcommand{\Gam}{\Gamma}
\newcommand{\Lam}{\Lambda}
\newcommand{\lam}{\lambda}
\newcommand{\om}{\omega}
\newcommand{\Om}{\Omega}
\newcommand{\ovG}{\overline G}
\newcommand{\si}{\sigma}


\def\mA{\mathbb{A}}
\def\mC{\mathbb{C}}
\def\mN{\mathbb{N}}
\def\mQ{\mathbb{Q}}
\def\mR{\mathbb{R}}
\def\mZ{\mathbb{Z}}
\def\mF{\mathbb{F}}
\def\cA{\mathcal{A}}
\def\cC{\mathcal{C}}
\def\cD{\mathcal{D}}
\def\cE{\mathcal{E}}
\def\cF{\mathcal{F}}
\def\cG{\mathcal{G}}
\def\cM{\mathcal{M}}
\def\cP{\mathcal{P}}
\def\cS{\mathcal{S}}
\def\cO{\mathcal{O}}
\def\cL{\mathcal{L}}
\def\cQ{\mathcal{Q}}
\def\cK{\mathcal{K}}
\def\cH{\mathcal{H}}
\def\cX{\mathcal{X}}
\def\i{\rm i\/}

\begin{abstract}
Using suitable deformations of simplicial trees and the duality
theory for median sets, we show that every free action on a median 
set can be extended to a free and transitive one. We also prove that
the category of median groups is a reflective full subcategory of
the category of free actions on pointed median sets.
\smallskip


\end{abstract}

\section*{Introduction}

$\quad$ Given an action $G \times \X \lra \X$ of a group $G$ on
a mathematical structure $\X$, a natural question to ask is whether it 
could be extended to a transitive action $\wh{G} \times \wh{\X} \lra
\wh{\X}$ on a structure $\wh{\X}$ of the same type with 
$\X$. Though in the simplest case when $\X$ is a set, this question
has a positive answer, the problem is not at all easy in the case
of actions on certain specific structures. In the present paper, 
having three parts, we consider the question above
in the frame of actions on generalized trees. 

Using the connection between Lyndon length functions and actions 
on $\Lam$-trees, as well as string-rewriting systems techniques, 
Ian Chiswell and Thomas M\" uller obtained quite recently the 
following nice result
\smallskip

{\bf Theorem} (cf. \cite[Theorem 5.4.]{Thomas}) {\em Let $G$ be a 
group acting freely and without inversions on a $\Lambda$-tree $\X$,
where $\Lam$ is a totally ordered abelian group.
Then there exists a group $\wh{G}$ acting freely, without inversions, 
and transitively on a $\Lambda$-tree $\wh{\X}$, together with a 
group embedding $G \to \wh{G}$ and a $G$-equivariant isometry} 
$\X \to \wh{\X}$.  
\smallskip

The main goal of the present work is to extend the result above in 
two directions: on the one hand, the $\Lambda$-trees, where $\Lam$ is
a totally ordered abelian group, are replaced by 
more general arboreal structures (median sets, in particular, 
faithfully full $\Lambda$-metric median sets (cf. \cite[1.3.]{RS}), where
$\Lambda$ is an abelian $l$-group); on the other hand, the free actions
without inversions are replaced by actions of a more general type
involving bicrossed products of groups. The method of proof is also different 
from that used in \cite{Thomas}: it is based on a suitable procedure of
deformation of simplicial trees into more general arboreal 
structures \footnote{$^)$ Compare with the procedure used in 
\cite[Proposition 4.2.]{DHAG} to extend a given locally linear
median group structure on a component of a suitable tree of groups
to the associated fundamental group.} as well as on the duality theory
for median sets (cf. \cite{Dual}).

The paper has three parts. The part I is devoted to free actions
on median sets. The results of part I will be applied in the part II 
to free actions on more general $\Lam$-trees, 
where $\Lam$ is an arbitrary abelian $l$-group, not necessarily
totally ordered, while the part III 
will be devoted to more general actions, not necessarily free, 
on generalized trees.

In order to state the main results of the part I of the paper,
we introduce (recall) some basic notions. More details will be given
in Section 1 having a preliminary character.
\smallskip

{\bf Definition 1.} By a {\em median set} (or {\em median algebra} 
\cite{Bandelt}), we understand a set $X$ together with a ternary 
operation $m : X^3 \lra X$ satisfying the following three equational axioms
\smallskip

{\bf (M 1) Symmetry}: $m(x, y, z) = m(y, x, z) = m(x, z, y)$

{\bf (M 2) Absorptive law}: $m(x, y, x) = x$

{\bf (M 3) Selfdistributive law}: $m(m(x, y, z), u, v) = 
m(m(x, u, v), m(y, u, v), z)$
\smallskip

The element $m(x, y, z)$ is called the {\em median of the triple} $(x, y, z)$.

Notice that {\bf (M 3)} can be replaced by

{\bf (M 3')} $m(m(x, u, v), m(y, u, v), x) = m(x, u, v)$
\smallskip

In particular, taking $u = y, v = z$, we obtain 
$m(x, y, m(x, y, z)) = m(x, y, z)$ for all $x, y, z \in X$.
For $x, y \in X$, we denote
$$\texttt{[} x, y \texttt{]} := \{z \in X\,|\,m(x, y, z) = z\} = 
\{m(x, y, z)\,|\,z \in X\}.$$
Notice that $\texttt{[} x, y \texttt{]} = \texttt{[} y, x \texttt{]}$,
and $u, v \in \texttt{[} x, y \texttt{]} \Lra \texttt{[} u, v \texttt{]} \sse
\texttt{[} x, y \texttt{]}$.
\smallskip

{\bf Definition 2.} Let $\X = (X, m)$ be a median set.

(1) A subset $M \sse X$ is said to be {\em convex} if 
$\texttt{[} x, y \texttt{]} \sse M$ for all $x, y \in M$. 
In particular, $\texttt{[} x, y \texttt{]}$ is convex
for all $x, y \in X$.

(2) A nonempty convex subset $M \sse X$ is {\em retractible} if for 
all $x \in X$ there exists (uniquely) $\psi(x) \in M$ such that 
$\texttt{[} x, a \texttt{]} \cap M = \texttt{[} \psi(x), a \texttt{]}$ 
for all (for some) $a \in M$; call the map $\psi : X \lra X$ with
$\psi(X) = M$ the {\em folding} associated to $M$.

(3) $\X$ is said to be {\em locally linear} if 
$\texttt{[} x, y \texttt{]} = \texttt{[} x, z \texttt{]} \cup 
\texttt{[} z, y \texttt{]}$ for all $x, y \in X, z \in \texttt{[} x, y \texttt{]}$.

(4) $\X$ is called {\em simplicial} (or {\em discrete}) 
if $\texttt{[} x, y \texttt{]}$ is finite for 
all $x, y \in X$.
\smallskip

Notice that, though median sets could seem too general to be considered
genuinely ``treelike'', they are, however, very well suited for the
study of various natural arboreal and metric structures on groups and rings -
see for instance \cite{DHAG}, \cite{RS} - \cite{arxiv}, 
\cite{Drutu2}, \cite{Drutu}. In particular, a key 
role is played by median groups \footnote{$^)$ The median sets and the 
median groups form proper subclasses of the larger class of
{\em connected median groupoids of groups} introduced in \cite{UC1, UC2}
as a frame for an extension of the Bass-Serre theory 
\cite{Trees} to actions on more general arboreal structures.}, defined as follows.
\smallskip

{\bf Definition 3.} By a {\em median group} \cite[2.2.]{AC} 
we understand a group $G$ together with a ternary operation 
$m : G^3 \lra G$ making $G$ a median set on which the 
group $G$ acts freely and transitively from the left, i.e.
$$\forall u, x, y, z \in G,\,u\,m(x, y, z) = m(u x, u y, u z)$$

The median group $(G, m)$ is called {\em locally linear} ({\em simplicial}) if
its underlying median set is locally linear (simplicial).
\smallskip

We denote by {\bf AMS} the category of actions on median sets. The
objects of {\bf AMS} are systems $(G, \X)$, where $G$ is a group,
together with an action from the left 
$G \times \X \lra \X$, $(g, x) \mapsto g \cdot x$ on a nonempty 
median set $\X = (X, m)$, while we take as morphisms 
$(G, \X) \lra (G', \X')$ the pairs $(\psi_0, \psi)$, where $\psi_0 : G \lra G'$
is a group morphism, and $\psi : \X \lra \X'$ is a morphism of 
median sets, compatible with the actions, i.e. $\psi(g \cdot x) = \psi_0(g) \cdot
\psi(x)$ for all $g \in G, x \in X$. We denote by {\bf FAMS} the full
subacategory of {\bf AMS} whose objects are the free actions on nonempty
median sets, while by {\bf FTAMS} we denote the full subcategory of
{\bf FAMS} consisting of the free and transitive actions on nonempty
median sets. 

We consider also the category {\bf FAPMS} of free actions on pointed
median sets having  as objects the systems $(G, \X, x_0)$, where 
$G$ is a group acting freely on a median set $\X = (X, m)$, while 
$x_0 \in X$ is a fixed base point.
The morphisms in {\bf FAPMS} are the morphisms in {\bf FAMS} which
preserve the base points. The category {\bf MG} of median groups, 
with naturally defined morphisms, is equivalent with {\bf FTAMS} and
is identified with a full subcategory of {\bf FAPMS}, by taking
the neutral element $1$ of a median group $\G = (G, m)$ as the base point
of its underlying median set.

A strong connection between free actions on median sets and median
groups is described by the first main result of the part I of the paper.
\smallskip

{\bf Theorem 1.} {\em Let $H$ be a group acting freely on a nonempty set 
$X$. Fix a basepoint $b_1 \in X$. We denote by $\cM(X)$ the set
of all median operations $m : X^3 \lra X$ which are compatible with
the action of $H$, i.e. $m(h x, h y, h z) = h \cdot m(x, y, z)$ for
all $h \in H, x, y, z \in X$. Let $\cM_l(X) (\cM_s(X))$ be the 
subset of $\cM(X)$ consisting of those median operations which
are locally linear (simplicial). Then there exists a group 
$\wh{H}$ containing $H$ as its subgroup, together with an embedding 
$\iota : X \lra \wh{H}$ and a retract $\varphi : \wh{H} \lra X$ 
such that the following hold.
\ben
\item[\rm (1)] The maps $\iota$ and $\varphi$ are $H$-equivariant,
i.e. $\iota(h \cdot x) = h \iota(x), \varphi(h u) = h \cdot \varphi(u)$
for all $h \in H, x \in X, u \in \wh{H}$. 
\item[\rm (2)] $\iota(b_1) = 1$, so $\iota(H b_1) = H$, and $\iota(X)$
generates the group $\wh{H}$.
\item[\rm (3)] Every $m \in \cM(X)$ extends uniquely to a ternary operation 
$\wh{m} : \wh{H}^3 \lra \wh{H}$ such that $(\wh{H}, \wh{m})$ is a
median group and the map $\iota \circ \varphi : \wh{H} \lra \wh{H}$ 
is a folding identifying $X$ with a retractible convex subset of 
the median set $(\wh{H}, \wh{m})$.
\item[\rm (4)] The median group $(\wh{H}, \wh{m})$ is locally linear
(simplicial) provided $m \in \cM_l(X)$ ($m \in \cM_s(X)$).
\een} 
\smallskip

The other two main results of the part I of the paper
are concerned with the construction of two functors
on the category {\bf FAMS} with suitable universal properties.
\smallskip

{\bf Theorem 2.} {\em Let $H$ be a group acting freely on a
nonempty median set $\X = (X, m)$. Then there exists a group
$\wh{H}$ acting freely on a median set $\wh{\X} = (\wh{X}, \wh{m})$,
together with a group {\em embedding} $\iota_0 : H \lra \wh{H}$ and
a $H$-equivariant {\em embedding} of median sets \\$\iota : \X \lra \wh{\X}$
such that $\iota(X) \sse \wh{H} \iota(x)$ for some (for all)
$x \in X$, and the following universal property is satisfied.

$({\rm RTUP})$ For every group $\wt{H}$ acting freely on a
median set $\wt{\X} = (\wt{X}, \wt{m})$ and for every
morphism $(\psi_0, \psi) : (H, \X) \lra (\wt{H}, \wt{\X})$ in
{\bf FAMS} such that $\psi(X) \sse \wt{H} \psi(x)$ for some (for all) 
$x \in X$, there exists uniquely a morphism 
$(\wh{\psi}_0, \wh{\psi}) : (\wh{H}, \wh{\X}) \lra
(\wt{H}, \wt{\X})$ in {\bf FAMS} such that} 
$\wh{\psi}_0 \circ \iota_0 = \psi_0, \wh{\psi} \circ \iota = \psi$.
\smallskip

Call the (unique up to a unique isomorphism) free action 
$(\wh{H}, \wh{\X})$, extending $(H, \X)$ and satisfying (RTUP), the
{\em relatively-transitive closure} of the free action $(H, \X)$.
\smallskip

By iterating the functorial construction above we obtain
\smallskip

{\bf Theorem 3.} {\em Let $H$ be a group acting freely on a
nonempty median set $\X$. Then there exists a
group $\fH$ acting {\em freely and transitively} on a median set $\fX$,
together with an {\em embedding} $(\iota_0, \iota) : (H, X) \lra (\fH, \fX)$
in the category {\bf FAMS} such that the following
universal property is satisfied.

$({\rm TUP})$ For every group $\wt{H}$ acting {\em freely and 
transitively} on a median set $\wt{\X}$ and for every morphism
$(\psi_0, \psi) : (H, \X) \lra (\wt{H}, \wt{\X})$ in the category
{\bf FAMS}, there exists uniquely a morphism
$(\Psi_0, \Psi) : (\fH, \fX) \lra (\wt{H}, \wt{\X})$ in {\bf FTAMS}
such that} $\Psi_0 \circ \iota_o = \psi_0, \Psi \circ \iota = \psi$.
\smallskip

Call the (unique up to a unique isomorphism) free and transitive 
action $(\fH, \fX)$ extending $(H, \X)$, with the property (TUP), the
{\em transitive closure} of the free action $(H, \X)$.
\smallskip

In category and median group theoretic terms, Theorem 3 can be 
rephrased as follows.
\smallskip

{\bf Theorem 3'.} {\bf MG} {\em is a reflective full subcategory
of {\bf FAPMS}, i.e. the embedding functor 
{\bf MG} $\lra$ {\bf FAPMS} has a left adjoint which {\em embedds}
every free action on a pointed median set into the median group
{\em freely generated} by it}.
\smallskip

The part I of the paper is organized in five sections. Section 1
introduces the reader to the basic notions and facts concerning 
median sets, free actions on median sets and median groups. In 
Section 2 we associate to an arbitrary free action of a group $H$
on a nonempty set $X$ a group $\wh{H}$ with an
underlying tree structure, together with a natural
embedding of the $H$-set $X$ into $\wh{H}$. Section 3
is devoted to the proof of a more explicit version of Theorem 1 by providing
a procedure for deformation of the underlying simplicial
tree of $\wh{H}$ introduced in Section 2 into suitable
median group operations on $\wh{H}$ which extend given median 
operations on the $H$-set $X$. Then, in Section 4 we use Theorem 1 
and the duality theory for median sets to prove Theorem 2.
Finally, by iteration of the functorial construction provided
by Theorem 2, we prove Theorem 3 in Section 5.


\section{Preliminaries on median sets, free actions on median sets
and median groups}

\subsection{Median sets}

$\quad$ The notion of median set (algebra) (see Definition 1 from
Introduction) appeared as a common generalization of trees and
distributive lattices. We recall here some basic definitions and
properties related to median sets. For proofs and further details 
we refer the reader to the papers \cite{Bandelt}, \cite{DF}, 
\cite{Dual}, \cite{UC1}, \cite{Birkhoff}, \cite{Isbell}, \cite{Roller}, 
\cite{Sholander}, \cite{Sholander2}.

The median sets form a category with naturally defined morphisms.

Let $\X = (X, m)$ be a nonempty median set. To any element $a \in X$ we
associate the binary operation 
$(x, y) \mapsto x \mathop {\cap }\limits_{a} y := m(x, a, y)$. 
With respect to this operation, $X$ is a meet-semilattice
with the induced partial order $x \mathop {\subset}\limits_{a} y \Llra
m(x, a, y) = x$ and the least element $a$, while for all $x, y, z \in X,
m(x, y, z)$ is the join with respect to $\mathop {\subset}\limits_{a}$
of the triple $(x \mathop {\cap }\limits_{a} y, y \mathop {\cap }\limits_{a} z, 
z \mathop {\cap }\limits_{a} x)$.

\begin{de} \em  
\ben
\item[\rm (1)] By a {\em median subset} of $\X$ we understand any 
subset $A$ of $X$ which is closed under the median operation $m$. 
\item[\rm (2)] A subset $A$ of $X$ is called {\em convex} if the 
following stronger condition is satisfied: for all $x, y, z \in X, 
m(x, y, z) \in A$ whenever $x, y \in A$.
\een
\end{de}

The intersection of any family of convex subsets of the median set 
$\X$ is convex too. Thus we may speak on the {\em convex closure of} 
(or {\em the convex subset generated by}) a subset $A \sse X$, and 
denote it by $\texttt{[} A \texttt{]}$. Notice that $\texttt{[} \es \texttt{]} = \es, 
\texttt{[} \{ a\}\texttt{]} = \{a \}$ for $a \in X$, and 
$\texttt{[} \{a, b \} \texttt{]} = \texttt{[} a, b \texttt{]}: =
\{ m(a, b, x)\,|\,x \in X \}= \{x \in X\,|\,m(a, b, x) = x\}$ for 
$a, b \in X$.

\begin{de} \em 
\ben
\item[\rm (1)] By a {\em cell} of $\X$ we mean a 
convex subset of the form $\texttt{[} a, b \texttt{]}$ with $a, b \in X$. 
\item[\rm (2)] Given a cell $C$ of $\X$, any 
$a \in X$ for which there exists $b \in X$ such that 
$C = \texttt{[} a, b \texttt{]}$ is called an {\em end} of the cell $C$. 
The (nonempty) subset of all ends of a cell $C$ is denoted by 
$\partial C$ and is called the {\em boundary} of $C$.
\een
\end{de}

The boundary $\partial C$ of a cell $C$ is 
a median subset of $C$, and there is a canonical automorphism 
$\neg$ of the median set $\partial C$ such that $C= \texttt{[} a, \neg \,a \texttt{]}$ 
and $\neg \neg \ a = a$ for all $a \in \partial C$. For any element 
$ a \in \partial C$, the cell $C$ becomes a bounded distributive lattice 
with the partial order $\mathop {\subset}\limits_{a}$, 
the meet $\mathop {\cap}\limits_{a}$, the join 
$\mathop {\cap}\limits_{\neg \ a}$, the least element $a$, and the 
last element $\neg \ a$. The boundary 
$\partial C$ is identified with the boolean subalgebra of the 
distributive lattice $(C, \mathop {\subset}\limits_{a})$ consisting 
of those elements which have (unique) complements. 

\begin{de} \em 
A median set $\X$ is said to be {\it locally boolean} if $C = \partial C$
for all cells $C$ of $\X$.
\end{de}

In particular, any boolean algebra $(L; \wedge, \vee, \neg, 0, 1)$
is a locally boolean median set with respect to the canonical median
operation
$$m(x, y, z) := (x \wedge y) \vee (y \wedge z) \vee (z \wedge x) =
(x \vee y) \wedge (y \vee z) \wedge (z \vee x).$$ 
For $x, y, z \in L$, $\texttt{[} x, y \texttt{]} = \texttt{[} x \wedge y,
x \vee y \texttt{]} = \texttt{[} m(x, y, z), m(x, y, \neg z) \texttt{]}$, 
in particular, $L = \texttt{[} 0, 1 \texttt{]} = \texttt{[} x, \neg x 
\texttt{]}$ for all $x \in L$.

\begin{rem} \em 
Let $\X = (X, m)$ be a median set. Then the following assertions hold.
\ben
\item[\rm (1)]
For $a, b, c \in X, \texttt{[} a, b \texttt{]} \cap \texttt{[} a, c \texttt{]} = 
\texttt{[} a, \mathop {b \cap c}\limits_{a}], \texttt{[} a, b \texttt{]} 
\cap \texttt{[} b, c \texttt{]} \cap \texttt{[} c, a \texttt{]} = \{ m(a, b, c)\}$, 
and $c \in \texttt{[} a, b \texttt{]} \Llra \texttt{[} a, c \texttt{]} \cap 
\texttt{[} b, c \texttt{]} = \{ c \}$. 
\item[\rm (2)] For nonempty finite subsets $A, B \sse X$,  
$\texttt{[} A \texttt{]} \cap \texttt{[} B \texttt{]}$ is 
either empty or the finitely generated convex subset 
$\texttt{[} \{ \mathop{\cap}\limits_{a} B\,|\, a \in A \} \texttt{]} = 
\texttt{[} \{ \mathop {\cap}\limits_{b} A\,|\, b \in B \} \texttt{]}$, 
where $\mathop {\cap}\limits_{a} B$ 
denotes the meet of the finite set $B$ with respect to
the partial order $\mathop {\subset}\limits_{a}$. In particular, for 
$a, b, c, d \in X, \texttt{[} a, b \texttt{]} \cap \texttt{[} c, d \texttt{]}$ 
is either empty or the cell $\texttt{[} a \mathop {\cap}\limits_{c} b, a \mathop 
{\cap}\limits_{d} b\texttt{]} = \texttt[c \mathop {\cap}\limits_{a} d, c \mathop 
{\cap}\limits_{b} d \texttt{]}$.
\een
\end{rem}

Among the convex subsets of a median set, we distinguish the 
prime ones, defined as follows.

\begin{de} \em
A convex subset $P$ of a median set $\X$ is said to be {\it prime} if 
its complement $X \sm P$ is also a convex subset of $\X$.
\end{de}

Thus the set $Spec\, \X$ of prime convex subsets of $\X$ is closed 
under the involution $P \mapsto \neg P := X \sm P$, and contains the empty 
set $\es$ as well as the whole $X$. For any subset $A$ of $X$, set 
$U(A) := \{P \in Spec\, \X\,|\, P \cap A = \es \}$,
$V(A) := \{ P \in  Spec\, \X\,|\,A \sse P \}$, and $U(a) := U(\{a\}),
V(a) := V(\{a\})$ for $a \in A$, whence $U(A) = \bigcap_{a \in A} U(a),
V(A) = \bigcap_{a \in A} V(a)$. The next result collects some basic
properties of the space $Spec\,\X$. For a proof see \cite[Theorems 5.2.1,
6.4.]{Dual}

\begin{te}
\ben
\item[\rm (1)] For $A, B \sse X, V(A) \cap U(B) \neq \es \Llra
[A] \cap [B] = \es$; in particular, 
$[A] = \mathop {\bigcap}\limits_{P \in V(A)} P$.
\item[\rm (2)] With respect to the topology defined by the basic
open sets $U(A)$ for $A$ ranging over the finite subsets of $X$,
$Spec\,\X$ is an irreducible spectral space \footnote{$^)$ 
A topological space $S$ is said to be {\em spectral} (or 
{\em coherent}) if

i) $S$ is {\em sober}, i.e. every irreducible nonempty closed 
subset of $S$ is the  closure of a unique point of $S$, and 

ii) the family of all quasicompact open subsets of $S$ is 
closed under finite intersection (in particular, $S$ itself is quasicompact) 
and forms a base for the topology of $S$.

A map $f : S' \lra S$ between spectral spaces is called {\em coherent} if
$f^{- 1}(U) \sse S'$ is a quasicompact open set provided $U \sse S$ is
a quasicompact open set. In particular, a coherent map is continuous.} 
with the generic point $\es$ and the unique closed point $X$.
\item[\rm (3)] The proper quasicompact open subsets of $Spec\,\X$ 
form a distributive lattice
$$\cL(\X) := \{\bigcup_{i = 1}^n  U(A_i)\,|\,n \geq 1, \es \neq A_i 
\sse X\,{\rm finite,\, for}\,i = \overline{1,\,n}\,\},$$ 
closed under the {\em negation} operator 
$$\cU \mapsto \neg\, \cU := \{P \in Spec\,\X)\,|\,\neg P \not \in \cU\},$$
while the embedding $X \lra \cL(\X), x \mapsto U(x)$ identifies 
the median set $\X$ to the invariant subset 
$\{\cU \in \cL(\X)\,|\,\neg\, \cU = \cU\}$, with the canonical
median operation
$$\fm(\cU_1, \cU_2, \cU_3) := (\cU_1 \cap \cU_2) \cup (\cU_2 \cap \cU_3) \cup (\cU_3 \cap \cU_1) =
(\cU_1 \cup \cU_2) \cap (\cU_2 \cup \cU_3) \cap (\cU_3 \cup \cU_1).$$
\item[\rm (4)] The correspondence $\X \mapsto Spec\, \X$ yields a 
duality between the category of median sets and a category of irreducible 
spectral spaces with a suitable additional structure \footnote{$^)$ More
precisely, according to \cite{Dual}, the objects of the dual category
of the category of median sets are the systems $(S, 0, 1, \neg)$, where
$S$ is an irreducible spectral space with generic point $0$, $1$ is the unique
closed point of $S$, and $\neg : S \lra S$ is an involution satisfying 
the following conditions.

(1) For every quasicompact open set $U \sse S$, the set
$\neg U := \{s \in S\,|\,\neg s \notin U\}$ is quasicompact open.

(2) The quasicompact open sets $U \sse S$ satisfying $\neg U = U$
generate the topology of $S$.

It follows that $\neg 0 = 1$.

The morphisms $f : (S, 0, 1, \neg) \lra (S', 0', 1', \neg')$ are
the coherent maps $f : S \lra S'$ satisfying $f(0) = 0'$
and $f \circ \neg = \neg' \circ f$, whence $f(1) = 1'$.}.

\item[\rm (5)] The correspondence $\X \mapsto \cL(\X)$ yields an
equivalence between the category of median sets and the category
of distributive lattices with negation $(L, \neg)$ which are generated 
(as lattices) by their invariant subsets $\{a \in L\,|\,\neg\,a = a\}$.
\een
\end{te}

As a corollary we obtain a description of the median set fms$(A)$ 
{\em freely generated by} an arbitrary set $A$.

\begin{co}
\ben
\item[\rm (1)] The restriction map $Spec\,{\rm fms}(A) \lra \cP(A),
P \mapsto P \cap A$ is bijective, identifying the spectral space
$Spec\,{\rm fms}(A)$ to the power set $\cP(A)$ with the basic open sets
$U(F) = \cP(A \sm F)$ for $F$ ranging over the finite subsets of $A$.
\item[\rm {2}] The elements of the distributive lattice $\cL({\rm fms}(A))$
correspond bijectively to families $(F_i)_{i = \overline{1,\,n}},\, n \geq 1$,
where the $F_i$'s are incomparable nonempty finite subsets of $A$; such a
family $(F_i)_{i = \overline{1,\,n}}$ corresponds to the proper quasicompact
open set $\dbigcup_{i = 1}^n  U(F_i) = \bigcup_{i = 1}^n \cP(A \sm F_i)$.
\item[\rm (3)] The negation operator sends a family $(F_i)_{i = \overline{1,\,n}}$
as above to the finite family of the subsets $E \sse \dbigcup_{i = 1}^n F_i$
which are minimal with the property $E \cap F_i \neq \es$ for $i = \overline{1,\,n}$.
\item[\rm (4)] The elements of the median set ${\rm fms}(A)$ freely generated by
the set $A$ correspond bijectively to families $(F_i)_{i = \overline{1,\,n}},\, n \geq 1$,
where the $F_i$'s are incomparable nonempty finite subsets of $A$ satisfying
\ben
\item[\rm (i)] $F_i \cap F_j \neq \es$ for $1 \leq i, j \leq n$, and
\item[\rm (ii)] for each subset $E \sse \dbigcup_{i = 1}^n F_i$
such that $E \cap F_i \neq \es$ for $i = \overline{1,\,n}$, there is
$1 \leq j \leq n$ such that $F_j \sse E$.
\een
\een
\end{co}

In particular, the finitely generated median sets are finite, and
any median set is a direct limit of finite median sets.

\begin{de} \em
A median set $\X$ is said to be {\em locally linear} if the 
following equivalent conditions are satisfied.
\ben
\item[\rm (1)] Every cell $C$ of $\X$ has at most two ends, 
i.e. $|\partial C| \in \{1, 2\}$.
\item[\rm (2)] For all $a, b \in X$, the cell $\texttt{[} a, b \texttt{]}$ 
is the set-theoretic union of the cells $\texttt{[} a, c \texttt{]}$ 
and $\texttt{[} b, c \texttt{]}$ provided $c \in \texttt{[} a, b \texttt{]}$.
\item[\rm (3)] For $a, b \in X$, the partial order 
$\mathop {\subset}\limits_{a}$ restricted to the cell 
$\texttt{[} a, b \texttt{]}$ is total (linear) with the 
least element $a$ and the last element $b$.
\item[\rm (4)] For $P, Q \in Spec\, \X$ such that $P \cap Q \neq \emptyset$ 
and $P \cup Q \neq X$, either $P \sse Q$ or $Q \sse P$.
\een
\end{de}

\begin{rem} \em
The locally linear median sets are strongly related with
order-trees. Recall that an {\em order-tree} is a poset $\T = (T, \leq)$ 
satisfying

(i) For every pair $(x, y)$ of elements in $T$, there exists the 
meet $x \wedge y$, and

(ii) $x \leq z$ and $y \leq z$ imply either $x \leq y$ or $y \leq x$.

Since for any triple $(x, y, z)$ of elements in an order-tree $\T$,
the set $\{x \wedge y, y \wedge z, z \wedge x\}$ has at most two distinct elements,
it follows that $\T$ has a natural structure of locally linear median 
set with $$m(x, y, z) := (x \wedge y) \vee (y \wedge z) \vee (z \wedge x) \in
\{x \wedge y, y \wedge z, z \wedge x\},$$
and $\texttt{[} x, y \texttt{]} = \texttt{[} x \wedge y, x \texttt{]} 
\cup \texttt{[} x \wedge y, y \texttt{]}$.
\end{rem}

The next notion will be very useful in the present paper.

\begin{de} \em By a {\em folding} of a median set $\X = (X, m)$ 
we mean a map $\varphi : X \lra X$ satisfying 
$\varphi(m(x, y, z)) =  m(\varphi (x), y, \varphi (z))$ for all 
$x, y, z \in X$.
\end{de}

One checks easily that a map $\varphi : X \to X$ is a folding 
if and only if $\varphi$ is an idempotent endomorphism of the median set $\X$ and 
the image $\varphi(X)$ is a convex subset of $\X$. In addition, according 
to \cite[Proposition 7.3.]{DF}, the map
$\varphi \mapsto \varphi(X)$ maps bijectively the foldings of 
$\X$  onto the nonempty convex subsets $A$ of 
$\X$ satisfying the following equivalent conditions.

i) $A$ is {\em retractible}, i.e. there is a (unique) median set 
retract $ p : X \lra A$ of the median set embedding 
$A \lra X$;

ii) For some (for all) $a \in A, A \cap \texttt{[} a, x \texttt{]}$ 
is a cell for all $x \in X$;

iii) For all $x \in X$, the meet $\mathop{\cap}\limits_{x} A$ with 
respect to the partial order $\mathop{\subset}\limits_{x}$ exists and 
belongs to $A$.

In particular, to any nonempty finite subset $A$ of a median set
$\X = (X, m)$, we associate the folding $\varphi_A$ defined by
$\varphi_A(x) = \mathop {\cap }\limits_{x} A$, the meet with respect to
the partial order $\mathop{\subset}\limits_{x}$ of the finite set $A$,
whence $\dbigcap_{a \in A} [x, a] = [x, \varphi_A(x)]$ for all $x \in X$,
and $\varphi_A(X) = [A]$.

\begin{de}
A median set $\X = (X, m)$ is said to be {\em simplicial} (or {\em discrete})
if every cell of $\X$ has finitely many elements.
\end{de}

To a simplicial median set $\X$ one assigns an integer-valued 
"distance" function $d : X \times X \lra \N$, 
where for $x, y \in X, d(x, y)$ is the length of some 
(of any) maximal chain in the finite distributive lattice 
$(\texttt{[} x, y \texttt{]}, \mathop{\subset}\limits_{x})$. 
With respect to $d$, $\X$ becomes a $\Z$-metric space such 
that for all $x, y \in X, \texttt{[} x, y \texttt{]} = 
\{z \in X\,|\,d(x, z) + d(z, y) = d(x, y)\}$, 
and the map $\texttt{[} x, y \texttt{]} \lra [0, d(x, y)], 
z \mapsto d(x, z)$, induced by $d$, is onto. In particular, 
$d(x, y) = d(u, v)$ whenever $\texttt{[} x, y \texttt{]} = 
\texttt{[} u, v \texttt{]}$, so we may speak on the "diameter" 
of any cell of $\X$. Notice also that for $x, y, z \in X, 
d(x, y \mathop{\cap}\limits_{x} z) = 
\frac{1}{2} (d(x, y) + d(x, z) - d(y, z))$. 

An equivalent graph theoretic description of simplicial 
median sets is given by \cite[Proposition 7.3]{UC1}. 
In particular, the customary {\em simplicial trees} 
(i.e. acyclic connected graphs) are identified with 
those simplicial median sets $\X$ which are 
locally linear, i.e. for all $x, y \in X$, the map 
$\texttt{[} x, y \texttt{]} \lra [0, d(x, y)]$ induced by 
$d$ is bijective. 

Notice also that in any simplicial median set, the nonempty 
convex subsets are retractible, and the finitely generated
convex subsets are finite.

\subsection{Groups acting freely on median sets}

$\quad$ By a {\em tree-free group} we mean a group having an action on a 
$\Lam$-tree, for some totally ordered abelian group $\Lam$, which is free 
and without inversions. This means that every non-identity element acts as 
a hyperbolic isometry (see \cite[Chapter 3]{Chiswell}). As any $\Lam$-tree,
where $\Lam$ is a totally ordered abelian group, has an underlying structure
of locally linear median set, the tree-free
groups form a remarkable subclass of the larger class ${\bf MSFG}$
consisting of the groups having a free action on some nonempty median
set.

\begin{lem}
${\bf MSFG}$ is a quasivariety of groups, i.e. it is closed under isomorphisms,
subgroups, and reduced products, and it contains the trivial group 
${\bf 1}$ (i.e. {\bf MSFG} $\neq \es$).
\end{lem}

\bp
We have only to show that the class ${\bf MSFG}$ is closed under
reduced products.  Let $I$ be a nonempty set, $\cF$ a filter on $I$,
and $(G_i)_{i \in I}$ a family of members of ${\bf MSFG}$. For 
each $i \in I$, let $G_i \times X_i \lra X_i$ be a free action on a
nonempty median set $X_i$. Using the filter $\cF$, we define the
normal subgroup $N$ of the product $\prod_{i \in I} G_i$ and
the congruence $\equiv$ on the median set product $\prod_{i \in I} X_i$
by
$$N := \{(g_i)_{i \in I} \in \prod_{i \in I} G_i\,|\,\{i \in I\,|\,g_i = 1_i\} \in \cF\},$$ 
$$x \equiv y \Llra \{i \in I\,|\,x_i = y_i\} \in \cF,\, 
{\rm\,for}\,x = (x_i)_{i \in I}, y = (y_i)_{i \in I} \in \prod_{i \in I} X_i.$$
One checks easily that the canonical free action $(\prod_{i \in I} G_i) \times
(\prod_{i \in I} X_i) \lra \prod_{i \in I} X_i$ induces a free action
of the reduced product $G := (\prod_{i \in I} G_i)/N$ on the quotient
median set $X := (\prod_{i \in I} X_i)/\equiv$, as desired.
\ep

\begin{lem}
As a full subcategory of the category ${\bf G}$ of groups, ${\bf MSFG}$ 
is reflective. The reflector, the left adjoint of the embedding 
${\bf MSFG} \lra {\bf G}$, sends a group $G$ to its quotient $G/N$, 
where $N$ is the smallest normal subgroup such that $G/N$ belongs 
to ${\bf MSFG}$.
\end{lem}

\bp
Given a group $G$, we denote by $\cN$ the set of those normal subgroups
$U$ of $G$ for which the quotient $G/U$ belongs to ${\bf MSFG}$. The
set $\cN$ is nonempty since $G \in \cN$. Set $N := \dbigcap_{U \in \cN} U$,
the kernel of the canonical morphism $G \lra \dprod_{U \in \cN} G/U$.
As ${\bf MSFG}$ is closed under products and subgroups, it follows that
$G/N$ belongs to ${\bf MSFG}$, whence $N$ is the least member of the
poset $\cN$ with respect to inclusion. The required adjunction
property is immediate.
\ep

As a quasivariety, the class of groups ${\bf MSFG}$ is axiomatized by
quasi-identities (see \cite[Theorem 2.25]{Burris}, \cite[Theorem 9.4.7]{Wilfrid}).
These quasi-identities turn out to be quite simple according to the next lemma.

\begin{lem}
Let $G$ be a group. Then the following assertions are equivalent.
\ben
\item[\rm (1)] $G$ belongs to ${\bf MSFG}$.
\item[\rm (2)] The canonical action of $G$ on the median set ${\rm fms}(G)$
freely generated by the underlying set of $G$ is free.
\item[\rm (3)] For all $g \in G$, either $g$ is of infinite order
or the order of $g$ is a power of $2$.
\een
\end{lem}

\bp
(1) $\Lra$ (3). Assume that $G$ acts freely on the nonempty median set $\X$,
and assume that there is $g \in G$ of prime order $p \neq 2$, say $p = 2 k + 1,
k \geq 1$. We have to get a contradiction. We may assume without loss that
$G$ is cyclic of order $p$, generated by $g$. Fix an element $x \in X$, so
$G x$, the $G$-orbit of $x$, has cardinality $p$. We denote by $\cF$ the 
(finite) set of all subsets $F \sse G x$ of cardinality $|F| = k + 1$. For any
$F \in \cF$, the set $U(F) = \{P \in Spec\,\X\,|\,P \cap F = \es\}$ is
a basic quasicompact open subset of the spectral space $Spec\,\X$ of
prime convex subsets of $\X$ (cf. Theorem 1.6.(2)). Consequently, the 
union $\cU := \dbigcup_{F \in \cF} U(F) = \{P \in Spec\,\X\,|\,|P \cap G x| \leq k\}
\in \cL(\X)$ is a proper quasicompact open subset of $Spec\,\X$. Moreover $\cU = \neg \cU =
\{P \in Spec\,\X\,|\,X \sm P \not \in \cU\}$, therefore, according to
Theorem 1.6.(3), there is an unique element $y \in X$ such
that $\cU = U(y) = \{P \in Spec\,\X\,|\,y \not \in P\}$.
Equivalently, by Theorem 1.6.(1), $y$ is the unique common element
of the convex closure $[F] = \dbigcap_{F \sse P \in Spec\,\X} P$ of $F$ in 
the median set $\X$ for
$F$ ranging over $\cF$. Since $g$ acts as a permutation on $\cF$, we deduce
that $g y = y$, contrary to the assumption that $G$ acts freely on $\X$.

(3) $\Lra$ (2). Assume that $G$ satisfies (3) and there is $g \in G \sm \{1\}$
such that $g x = x$ for some $x \in {\rm fms}(G)$. We have to 
obtain a contradiction. By Corollary 1.7.(4), the element 
$x = g x \in {\rm fms}(G)$ is uniquely determined by a suitable finite 
family $\cF = (F_i)_{i = \overline{1,\,n}},\, n \geq 1$, of incomparable 
nonempty finite subsets of $G$, whence $g F_i = \{g h\,|\,h \in F_i\} =
F_{\si(i)}, i = \overline{1,\,n}$, for some permutation $\si$ of the finite set 
$\{1, \dots, n\}$. Consequently, $g$ has finite order, so we may assume
without loss that $g$ has order $2$. Set $F := \dbigcup_{i = 1}^n F_i$.
As $g F = \{g h\,|\,h \in F\} = F$, there is $E \sse F$ such that $F$ is
the disjoint union of its subsets $E$ and $g E$, in particular, $|F| = 2 k$
with $k = |E| = |g E| \geq 1$. We distinguish the following two cases.

(i) $E \cap F_i = \es$ for some $1 \leq i \leq n$. Then $F_i \sse g E$,
whence $F_{\si(i)} = g F_i \sse g^2 E = E$, therefore $F_i \cap F_{\si(i)} =
\es$, contrary to the condition (i) of Corollary 1.7.(4) satisfied by the
family $\cF$.

(ii) $E \cap F_i \neq \es$ for all $i = \overline{1,\,n}$. Then, by condition
(ii) of Corollary 1.7.(4) satisfied by the family $\cF$, there is 
$1 \leq j \leq n$ such that $F_j \sse E$, whence $F_{\si(j)} = g F_j \sse g E$,
so $F_j \cap F_{\si(j)} = \es$, again a contradiction.

Finally notice that the implication (2) $\Lra$ (1) is trivial.
\ep

\begin{co}
The quasivariety of groups {\bf MSFG} is axiomatized by the quasi-identities
$$x^p = 1 \rightarrow x = 1$$
for $p$ ranging over the set of odd prime numbers.
\end{co}

\subsection{Median groups}

$\quad$ We recall here some basic definitions and properties
related to median groups. For proofs and further details we
refer the reader to the papers \cite{DHAG}, \cite{AC}, \cite{arxiv}.

\begin{de} \em 
\ben
\item[\rm (1)] Let $G$ be a group. By a {\em median group operation} on
$G$ we understand a ternary operation $m : G^3 \lra G$ satisfying 
\ben
\item[\rm (i)] $(G, m)$ is a median set, and
\item[\rm (ii)] $u\, m(x, y, z) = m(u x, u y, u z)$ for all
$u, x, y, z \in G$.
\een
\item[\rm (2)] By a {\em median group} we understand a group
$G$ together with a median group operation $m$ on $G$.
\item[\rm (3)] By a {\em formally-median group} we mean a group
$G$ satisfying the following equivalent conditions.
\ben
\item[\rm (i)] There exists a median group operation on $G$.
\item[\rm (ii)] $G$ acts freely and transitively on some 
nonempty median set.
\item[\rm (iii)] There exists a $G$-equivariant retract
$\varphi : {\rm fms}(G) \lra G$ of the canonical $G$-equivariant
embedding of $G$ into the median set $({\rm fms}(G), m)$ freely
generated by the set $G$ such that $\varphi(m(x, y, z)) =
\varphi(m(\varphi(x), y, z))$ for all $x, y, z \in {\rm fms}(G)$.
\een
\een
\end{de}

Let $\G = (G, m)$ be a median group. Taking the neutral element
$1$ as a basepoint of the underlying median set of $\G$, we 
get the meet-semilattice operation $x \cap y = m(x, 1, y)$
with the induced partial order $\subset$. Thus $1 \subset x$ for all 
$x \in G$, and $x \subset y \Llra x \in \texttt{[} 1, y \texttt{]}$.
Notice also that $z \in \texttt{[} x, y \texttt{]} \Llra x^{- 1} z
\subset x^{- 1} y$, and 
$$m(x, y, z) = x (x^{- 1} y \cap x^{- 1} z) = y (y^{- 1} x \cap 
y^{- 1} z) = z (z^{- 1} x \cap z^{- 1} y)$$
for all $x, y, z \in G$. In particular, for $x, y \in G$, 
$x \cap y$ is the unique element $z \in G$ satisfying 
$z \subset x, z \subset y$ and $x^{- 1} z \subset x^{- 1} y$.

The next statement furnishes an useful order theoretic
description of median groups.

\begin{pr} {\em (cf. \cite[Proposition 2.2.1.]{AC})}
Let $G$ be a group. Then the map sending a ternary operation
$m : G^3 \lra G$ to the binary operation $\cap$, defined by
$x \cap y := m(x, 1, y)$, maps bijectively the median group
operations on $G$ onto the binary operations $\cap$ on $G$ satisfying
\ben
\item[\rm (1)] $(G, \cap)$ is a meet-semilattice; let $x \subset y
\Llra x \cap y = x$ be the induced partial order,
\item[\rm (2)] $1 \subset x$ for all $x \in G$,
\item[\rm (3)] $x \subset y, y \subset z \Lra z^{- 1} y \subset
z^{- 1} x$, and
\item[\rm (4)] $x^{- 1} (x \cap y) \subset x^{- 1} y$ for all $x, y \in G$.
\een
\end{pr}

In both the signatures $(1, ^{-1}, \cdot, m)$ and 
$(1, ^{-1}, \cdot, \cap)$ of type $(0, 1, 2, 3)$ and $(0, 1, 2, 2)$ 
respectively, the median groups form a variety.
In particular, the class of median groups is closed under 
arbitrary products, with group and median operations defined
component-wise.

\begin{de} \em 
A median group is said to be {\em locally linear} ( {\em simplicial} )
if its underlying median set is locally linear ( simplicial ).
\end{de}

\begin{co} {\em (cf. \cite[Corollary 2.2.2.]{AC})} Let $G$ 
be a group. Then the map sending a ternary operation 
$m : G^3 \lra G$ to the binary relation 
$x \subset y \Llra m(x, 1, y) = x$ maps bijectively
the locally linear median group operations on $G$ onto the
partial orders $\subset$ on $G$ satisfying
\ben
\item[\rm (1)] $1 \subset x$ for all $x \in G$,
\item[\rm (2)] $x \subset y, y \subset z \Lra z^{- 1} y \subset
z^{- 1} x$,
\item[\rm (3)] for all $x, y \in G$ there exists $z \in G$ such 
that $z \subset x, z \subset y$ and $x^{- 1} z \subset x^{- 1} y$, and
\item[\rm (4)] $x \subset z, y \subset z \Lra$ either $x \subset y$
or $y \subset x$.
\een
\end{co}

In particular, if $(G, \leq)$ is a {\em left-ordered} group \footnote{$^)$
A {\em left-ordered group} is a group $G$, together with a total order
$\leq$ on $G$ such that $u \leq v \Lra g u \leq g v$ for all $g, u, v \in G$.
A group $G$ is {\em left-orderable} if $(G, \leq)$ is left-ordered for
some total order $\leq$ on $G$.}, then the total order $\leq$ determines
a locally linear median group operation $m$ on $G$ whose associated
partial order $\subset$ is given by $u \subset v \Llra$ either
$1 \leq u \leq v$ or $v \leq u \leq 1$. In other words, the median
operation $m$ is induced by the betweenness relation associated to
the total order $\leq$. Notice that there exist groups $G$ together with
total orders $\leq$ such that $(G, \leq)$ is not left-ordered, but
the median operation induced by the betweenness relation associated to
the total order $\leq$ is compatible with the left multiplication
(see Remark 1.21.(2)). In any case, possible connections with the
much more studied class of left-orderable groups could be fruitful.
 
\begin{co} The necessary and sufficient condition for a median group
 $\G = (G, m)$ to be simplicial is that for all $x \in G$, the
cell {\em \texttt{[}} $1, x $ {\em \texttt{]}} 
$ = \{y \in G\,|\,y \subset x\}$
has finitely many elements.
\end{co}

\begin{rem} \em According to Corollary 1.15., the quasivariety {\bf MSFG} of
groups acting freely on median sets is described by very simple quasi-identities.
By contrast, the proper subclass of {\bf MSFG} consisting of the formally-median
groups is not enough investigated. Some particular classes of such groups are 
studied in \cite{DHAG}, \cite{AC}, \cite{RT}, \cite{arxiv}. For convenience
of the reader we discuss here only the simplest case of formally-median
cyclic groups. 

(1) Though all cyclic groups of order $2^n, n \in \N$, belong to 
{\bf MSFG}, only three of them with $n = 0, 1, 2$ are formally-median groups.
The corresponding median group operations are uniquely determined : the 
point, the segment $\texttt{[} 1, \si \texttt{]}, \si^2 = 1$, and the 
square $\texttt{[} 1, \si^2 \texttt{]} = \texttt{[} \si, \si^3 \texttt{]},
\si^4 = 1$, respectively.

(2) The infinite cyclic group $(\Z, +)$ has a canonical structure of
simplicial and locally linear median group with respect to the 
median group operation $m_0$ associated to the usual simplicial tree 
on $\Z$, induced by the natural order : 
$m_0(x, y, z) = y \Llra$ either $x \leq y \leq z$ or 
$z \leq y \leq x$. However there are still two
distinct median group operations $m_1$ and $m_{- 1}$  on $(\Z, +)$,
which are both locally linear but not simplicial, related each to other 
through the unique nonidentical automorphism $n \mapsto - n$ of the group $(\Z. +)$ : 
$$m_{- 1}(x, y, z) = - m_1(- x, - y, - z)\,\, {\rm for}\,\, x, y, z \in \Z.$$
Thus there are only two (up to isomorphism) median groups with the
underlying group $(\Z, +)$. To prove the assertion above and describe
explicitely the median group operation $m_1$, introduced in 
\cite[Remarks 3.2.(3)]{arxiv}, we proceed as follows. Let $m$ be
a median group operation on $(\Z, +)$, with the associated meet-semilattice
operation $\cap$ and partial order $\subset$ (cf. Proposition 1.17.)

Let $x, y \in \Z$ be such that $0 < x < y$, and let $z := m(0, x, y) = x \cap y$.
Assuming that $y < z$ and using the fact that the translation $n \mapsto n + 1$
is an automorphism of the median set $(\Z, m)$, we deduce that
$\N$ is contained in the median subset of $(\Z, m)$ generated by
the finite set $\{0, 1, \dots, z\}$, i.e. a contradiction. 
Similarly, assuming that $z < 0$ and using the translation 
$n \mapsto n - 1$, we deduce that $\Z_{\leq y}$ is contained in
the median subset generated by the finite set $\{z, z + 1, \dots, y\}$,
again a contradiction. Thus the following implication holds 
$$(\ast) \quad 0 < x < y \Lra 0 \leq x \cap y \leq y.$$
In particular, we get $1 \cap 2 = m(0, 1, 2) \in \{0, 1, 2\}$. We
distinguish the following three cases.
\ben
\item[\rm (i)] $1 \cap 2 = 1$, i.e. $1 \subset 2$. To obtain $m = m_0$,
it suffices by Proposition 1.17. to show that $n \subset n + 1$ for all $n \geq 1$. As 
$1 \subset 2$ by hypothesis, assuming by induction that $k \subset k + 1$
for $1 \leq k < n, n \geq 2$, we have to show that $s := m(0, n, n+ 1) = n$.
By $(\ast)$, $0 \leq s \leq n + 1$. Assuming that $s \neq n$, there are
three possibilities.
\ben
\item[\rm (a)] $1 \leq s \leq n - 1$, whence $s - 1 \subset n - 1 \subset n$ by
the induction hypothesis. On the other hand, $s \in \texttt{[} n, n + 1 \texttt{]}$,
therefore $s - 1 \in \texttt{[} n - 1, n \texttt{]}$, so $n - 1 \subset s - 1 \subset n$.
Consequently, $s - 1 = n - 1$, contrary to the assumption $s \leq n  - 1$.
\item[\rm (b)] $s = 0$, i.e. $0 \in \texttt{[} n, n + 1 \texttt{]}$, and
hence $n - 1 \in \texttt{[} 0, n \texttt{]} \sse \texttt{[} n, n + 1 \texttt{]}$.
As $1 \in \texttt{[} 0, 2 \texttt{]} \Lra n \in \texttt{[} n- 1, n + 1 \texttt{]}$,
we deduce that $n - 1 = n$, i.e. a contradiction.
\item[\rm (c)] $s = n + 1$, i.e. $n + 1 \in \texttt{[} 0, n \texttt{]}$,
therefore $\texttt{[} n - 1, n + 1 \texttt{]} \sse \texttt{[} 0, n \texttt{]}$.
Setting $t := m(0, n - 1, n + 1)$, we obtain by Remark 1.4.(2)
$$\texttt{[} n - 1, n + 1 \texttt{]} = \texttt{[} n - 1, n + 1 \texttt{]} \cap
\texttt{[} 0, n \texttt{]} = \texttt{[} t, m(n, n - 1, n + 1) \texttt{]} =
\texttt{[} t, n \texttt{]},$$
whence $t \notin \{n - 1, n, n + 1\}$, and hence 
$0 \leq t \leq n - 2$ by $(\ast)$. Assuming that $t \neq 0$, it follows
by the induction hypothesis that $t - 1 \subset n - 2 \subset n - 1 \subset n$.
Since $\texttt{[} n - 1, n + 1 \texttt{]} = \texttt{[} t, n \texttt{]} \Lra
\texttt{[} n - 2, n \texttt{]} = \texttt{[} t - 1, n - 1 \texttt{]}$ (by
translation with $- 1$), we deduce that $n - 1 = n$, i.e. a contradiction.
Consequently, $t = 0$, and hence $1 \in \texttt{[} 0, n \texttt{]} = 
\texttt{[} n - 1, n + 1 \texttt{]}$, whence 
$0 \in \texttt{[} n - 2, n \texttt{]}$ (by translation
with $- 1$). As $n - 2 \subset n$, it follows that 
$n - 2 = 0$, i.e. $\texttt{[} 0, 2 \texttt{]} = \texttt{[} 1, 3 \texttt{]}$,
whence $\texttt{[} 0, 2 \texttt{]} = \texttt{[} 2, 4 \texttt{]}$ (by 
translation with $1$), therefore $0 = 4$, i.e. again a contradiction. 

Since in all three possible cases (a), (b), (c) we obtain a 
contradiction, we deduce that $n \subset n + 1$ for all $n \in \N$, 
and hence $m = m_0$ as desired.
\een
\item[\rm (ii)] $1 \cap 2 = 2$, i.e. $2 \in \texttt{[} 0, 1 \texttt{]}$,
and hence $1 \in \texttt{[} 0, - 1 \texttt{]}$ (by translation with $- 1$),
whence $2 \subset 1 \subset - 1$. Applying succesively the translation
$k \mapsto k + 2$, it follows that $2 n \subset 2 n + 2 \subset 2 n + 1 
\subset 2 n - 1$ for $n \in \N$. Consequently, the median group operation
$m := m_1$ is uniquely determined by the betweenness relation induced by
the total order $\prec$ (or its opposite) on $\Z$, 
defined by $x \prec y$ if and only if one of the following assertions
hold

(a) $x, y$ are even and $x \leq y$;

(b) $x, y$ are odd and $y \leq x$;

(c) $x$ is even and $y$ is odd.
\item[\rm (iii)] $1 \cap 2 = 0$, i.e. $0 \in \texttt{[} 1, 2 \texttt{]}$,
whence $- 1 \in \texttt{[} 0, 1 \texttt{]}$ (by translation with $- 1$),
i.e. $- 1 \subset 1$. It follows that the median group operation $m := m_{- 1}$
is the conjugate of $m_1$ by the group automorphism $n \mapsto - n$,
and hence it is uniquely determined by the betweenness relation induced 
by the total order $\prec'$ (or its opposite) on $\Z$ obtained from $\prec$ by 
replacing (c) with

(c') $x$ is odd and $y$ is even.
\een

Notice that though the total orders $\prec$ and $\prec'$
are not compatible with the group operation on $\Z$, the induced median
operations $m_1$ and $m_{- 1}$ are so.

Notice also that $2 \Z$ is a median subgroup of $(\Z, m_i), i = 0, 1, - 1$,
and $m_i|_{2 \Z} = m_0|_{2 \Z}, i = 1, - 1$. However, by contrast with
the median group $(\Z, m_0)$ which has no proper convex subgroups, 
and hence no proper quotients, $2 \Z$ is the unique
proper convex subgroup of $(\Z, m_1) \cong (\Z, m_{- 1})$, inducing
the surjective morphism of {\em median groups} $(\Z, m_1) \cong (\Z, m_{- 1})
\lra \Z/2$, whose kernel is isomorphic to the median group $(\Z, m_0)$. In
other words, the isomorphic median groups $(\Z, m_1)$ and $(\Z, m_{- 1})$
are extensions of the median group $\Z/2$ by the median group $(\Z, m_0)$.

(3) The construction above has a nice interpretation in terms of
the nonstandard arithmetic. Let $^{\ast}\Z$ be an {\em enlargement} of
$\Z$. For our purposes it suffices to take $^{\ast}\Z$ an ultrapower
of $\Z$ relative to a nonprincipal ultrafilter on $\N$. We denote
by $^{\ast}m_i$ the median group operation on $^{\ast}\Z$
which extends the median group operation $m_i, i = 0, 1, - 1$.
Let $T_i := \{t \in\, ^{\ast}\Z\,|\,\forall x, y \in \Z,\, ^{\ast}m_i(t, x, y) \in
\Z\}$. $T_i$ is the maximal median subset of 
$(^{\ast}\Z,\, ^{\ast}m_i)$ lying over $\Z$ with the property that
$\Z$ is convex in $(T_i,\, ^{\ast}m_i|_{T_i})$.
Since $(\Z, m_0)$ is simplicial, it follows that $T_0 =\, ^{\ast}\Z$,
while 
$$T_1 = \{- 2 t\,|\,t \in\, ^{\ast}\N \sm \N\} \bigsqcup\, \Z
\,\bigsqcup \{- 2 t + 1\,|\,t \in \,^{\ast}\N \sm \N\},$$
$$T_{- 1} = \{2 t + 1\,|\,t \in\,^{\ast}\N \sm \N\} \bigsqcup\, \Z
\,\bigsqcup \{2 t\,|\,t \in\, ^{\ast}\N \sm \N\}$$
are proper {\em median submonoids} of the median group 
$(^{\ast}\Z, +,\,^{\ast}m_i)$ for $i = 1, - 1$ respectively,
containing $\Z$ as the maximal (convex) subgroup. 

Let us consider the congruence $\equiv_i$ on 
the median set $(T_i,\, ^{\ast}m_i|_{T_i})$, 
defined by 
$$t \equiv_i t' \Llra \forall x, y \in \Z,
\,^{\ast}m_i(t, x, y) = \,^{\ast}m_i(t', x, y).$$
According to \cite[A.1.Proposition]{DF},
the quotient median set $T_i/\equiv_i$ is isomorphic to the median set
Dir$(\Z, m_i)$ of the {\em directions} on the median set 
$(\Z, m_i)$, containing $\Z$ as the convex subset of {\em internal
directions}. It follows that Dir$(\Z, m_i) = 
\texttt{[} D_i, D_i' \texttt{]} = \{D_i\} \bigsqcup \Z \bigsqcup \{D_i'\}$, 
where the {\em external directions} $D_i, D_i'$ are the equivalence 
classes $\{- t\,|\,t \in\,^{\ast}\N \sm \N\},
\{t\,|\,t \in\,^{\ast}\N \sm \N\}$ for $i = 0$, 
$\{- 2 t\,|\,t \in\, ^{\ast}\N \sm \N\},
\{- 2 t + 1\,|\,t \in \,^{\ast}\N \sm \N\}$ for $i = 1$, and
$\{2 t + 1\,|\,t \in\,^{\ast}\N \sm \N\}, \{2 t\,|\,t \in\, ^{\ast}\N \sm \N\}$
for $i = - 1$ respectively. The induced total order on $\Z$ from 
$D_i$ to $D_i'$ is $\leq$ for $i = 0$, $\prec$ for $i = 1$, 
and $\prec'$ for $i = - 1$. The canonical free action of $\Z$ on $T_i$,
$(n, t) \in \Z \times T_i \mapsto n + t \in T_i$ induces an action
on Dir$(\Z, m_i) \cong T_i/\equiv_i$ which is obviously free and transitive on 
the set $\Z$ of internal directions, identical on the external
directions $D_0, D_0'$, and acting as $\Z/2$ on the pair of 
external directions $(D_i, D_i')$ for $i = 1, - 1$.
\end{rem}

We end this preliminary section with an useful lemma relating
free actions on median sets and median groups.

\begin{lem}
Let $G$ be a group, $H$ a subgroup of $G$, and $X \sse G$ a set of 
generators of $G$ such that $H \sse X$, and $H X = X$, whence 
$H \times X \lra X,(h, x) \mapsto h x$ is a free action of $H$ 
on the nonempty set $X$, and the embedding $\iota : X \lra G$ is $H$-equivariant. 
Let $\varphi : G \lra X$ be a $H$-equivariant retract of $\iota$.

We denote by $\cM(X)$ the set of the median operations $m : X^3 \lra X$
which are compatible with the action of $H$, i.e. $m(h x, h y, h z) =
h m(x, y, z)$ for $h \in H, x, y, z \in X$.

On the other hand, we denote by $\cM(G, \varphi)$ the set of 
those median group operations $\wh{m} : G^3 \lra G$ for which the map 
$\varphi$ is a folding, so $X$ is a retractible convex subset
of $(G, \wh{m})$ with associated folding $\varphi$. 

Then the restriction map $\cM(G, \varphi) \lra \cM(X)$ is injective.
\end{lem}
 
\bp
Let $I \sse X$ be a system of representatives of the $H$-orbits.
Assume that $1 \in I$, and set $I' := I \sm \{1\}$. Thus the disjoint
union $(H \sm \{1\}) \sqcup I'$ generates the group $G$. Let
$\wh{m} \in \cM(G, \varphi)$, and $m \in \cM(X)$ be its restriction.
To prove that $\wh{m}$ is the unique prolongation of $m$, it suffices 
to show by duality (cf. Theorem 1.6.(4)) that $Spec(G, \wh{m})$ is
uniquely determined by $Spec(X, m)$ and $\varphi$. 

The $H$-equivariant morphisms of median sets $\iota : (X, m) \lra (G, \wh{m})$
and \\ $\varphi : (G, \wh{m}) \lra (X, m)$ satisfying $\varphi \circ \iota =
1_X$ induce by duality the morphisms of spectral spaces
$$Spec(G, \wh{m}) \lra Spec(X, m), P \mapsto P \cap X, 
Spec(X, m) \lra Spec(G, \wh{m}), \fp \mapsto \varphi^{- 1}(\fp)$$
such that $\varphi^{- 1}(\fp) \cap X = \fp$ for all $\fp \in Spec(X, m)$.

As $\wh{m}$ is a median group operation, $G$ acts from the right on
$Spec(G, \wh{m}), (P, g) \mapsto P^g := g^{- 1} P$, while $H$ acts 
from the right on $Spec(X, m), (\fp, h) \mapsto \fp^h := h^{- 1} \fp$,
and the induced morphisms of spectral spaces are $H$-equivariant. 

Let $\cS := \{\varphi^{- 1}(\fp)^g = g^{- 1} \varphi^{- 1}(\fp)\,|\,
\fp \in Spec(X, m), g \in G\}$. Notice that $\es, G \in \cS$, and
$P \in \cS \Lra G \sm P \in \cS$. The inclusion $\cS \sse Spec(G, \wh{m})$ 
is obvious, so it remains to prove the opposite inclusion. 
Let $P \in Spec(G, \wh{m}) \sm \{\es, G\}$. We
distinguish the following three cases.
\smallskip

(1) $\fp := P \cap X \neq \es, X$. Then $P = \varphi^{- 1}(\fp) \in \cS$.
Indeed, assuming the contrary, say $P \not \sse \varphi^{- 1}(\fp)$,
let $g \in P$ be such that $\varphi(g) \not \in \fp$. As $\fp \neq \es$
by assumption, choose some $x \in \fp \sse P$. Since $\varphi$ is a folding
of the median set $(G, \wh{m})$ with $\varphi(G) = X$, and $P$ is convex, 
we deduce that $\varphi(g) \in [x, g] \cap X \sse P \cap X = \fp$, 
i.e. a contradiction. The case $\varphi^{- 1}(\fp) \not \sse P$ 
follows similarly by replacing $P$ with $G \sm P$.
\smallskip

(2) $X \sse P$. As $P \neq G$, choose an element $g \in G \sm P$
of minimal length $l(g)$ over the alphabet $J := (H \sm \{1\}) \sqcup I'^{\pm 1}$.
In particular $g \neq 1$, i.e. $l(g) \geq 1$, since $1 \in X \sse P$ by assumption. Let
$g = g' t$ be a reduced word representing $g$, with $t \in J$. Notice 
that $g' = g t^{- 1} \in P$ since $l(g') < l(g)$.
There are two possibilities.
\smallskip

(i) $t \in (H \sm \{1\}) \sqcup I'^{- 1}$. Then $t^{- 1} \in g^{- 1} P \cap X$,
while $1 \in X \sm g^{- 1} P$. Consequently, 
$\fp := g^{- 1} P \cap X \in Spec(X, m) \sm \{\es, X\}$,
therefore $P = g \varphi^{- 1}(\fp) \in \cS$ by (1).
\smallskip

(ii) $t \in I'$. Then $t = g'^{- 1} g \in X \sm g'^{- 1} P$ and
$1 = g'^{- 1} g' \in X \cap g'^{- 1} P$, whence $\fp := X \cap g'^{- 1} P \in
Spec(X, m) \sm \{\es, X\}$, and hence $P = g' \varphi^{- 1}(\fp) \in \cS$ by
(1) again.
\smallskip

(3) $P \cap X = \es$. Then $G \sm P \in \cS$ by (2), whence
$P \in \cS$ as desired.
\ep

In particular, taking $H = 1$, we obtain

\begin{co}
Let $G$ be a group, $X \sse G$ a set of generators with $1 \in X$,
and $\varphi : G \lra X$ a surjective map such that $\varphi(x) = x$
for all $x \in X$.

We denote by $\cM(X)$ the set of median operations on $X$, and by
$\cM(G, \varphi)$ the set of median group operations on $G$ for which
$X$ is a retractible conves subset with associated folding $\varphi$.

Then the restriction map $\cM(G, \varphi) \lra \cM(X)$ is injective. 
\end{co}

With the notation from Corollary 1.23., call the surjective map
$\varphi : G \lra X$ {\em admissible} if $\cM(G, \varphi) \neq \es$.
Let us give some simple examples of admissible maps.

\begin{exs} \em 
\ben
\item[\rm (1)] For $G = \la g\,|\,g^4 = 1 \ra \cong \Z/4, X = \{1, g\}$,
the map $\varphi : G \lra X$, with $\varphi(g^2) = \varphi(g) = g$,
$\varphi(g^3) = \varphi(1) = 1$, is the unique admissible map, and
$\cM(G, \varphi)$ consists of the unique median group operation on $G$ - 
the square with the pairs of opposite vertices $(1, g^2)$ and $(g, g^3)$ -, 
so the restriction map $\cM(G, \varphi) \lra \cM(X)$ is obviously bijective.
\item[\rm (2)] For $G = (\Z, +), X = \{0, 1\}$, the surjective map $\varphi : G \lra X$
with $\varphi^{- 1}(1) = \Z_{\geq 1}$ is the unique admissible map, and
$\cM(G, \varphi) = \{m_0\}$, where $m_0$ is the canonical median group operation
on $\Z$, corresponding to the natural simplicial tree on $\Z$, so
the restriction map $\cM(G, \varphi) \lra \cM(X)$ is obviously bijective.
\item[\rm (3)] A more interesting case is
$G = (\Z, +), X = \N$, where we have to find those median group operations
$m$ on $\Z$ satisfying the strong condition that the submonoid $(\N, +)$
is a retractible convex subset of $(\Z, m)$, whence, by translation with
elements $n \in \Z$, $\Z_{\geq n}$ is also retractible convex in $(\Z, m)$.
This task is easy since we already know, according to Remark 1.21.(2), that
there are only three distinct median group operations $m_0, m_1, m_{- 1}$ on
$(\Z, +)$. We see that only two of them, namely $m_0$ and $m_1$, satisfy the
requirement, providing the admissible maps $\varphi_i : \Z \lra \N, i = 0, 1$,
defined by 
\begin{center}
$\varphi_0(n) = \left\{\begin{array}{lcl}
n & \mbox{if} & n \geq 0 \\
0 & \mbox{if} & n < 0,
\end{array}
\right.$
\end{center}
with $\cM(\Z, \varphi_0) = \{m_0\}$, and 
\begin{center}
$\varphi_1(n) = \left\{\begin{array}{lcl}
n & \mbox{if} & n \geq 0 \\
0 & \mbox{if} & n < 0\,{\rm and}\,n \in 2 \Z \\
1 & \mbox{if} & n < 0\,{\rm and}\,n \in 2 \Z + 1,
\end{array}
\right.$
\end{center}
with $\cM(\Z, \varphi_1) = \{m_1\}$. Notice that $\varphi_1$ is the folding 
$n \mapsto m_1(0, n, 1)$ associated to the linear cell 
$\texttt{[} 0, 1 \texttt{]}_{m_1} = \N$. Thus, with the exception of
$m_0|_\N$ and $m_1|_\N$, the infinitely many median operations on
the countable set $\N$ do not extend to median group operations on $(\Z, +)$.
To extend them to suitable median group operations we are forced to 
forget the monoid structure of $\N$ and look for larger groups 
containing the countable {\em set} $\N$ (see Corollary 3.6).
\een
\end{exs}


\section{Simplicial trees induced by free actions on sets}

$\quad$ The main goal of this section is to explore the
underlying simplicial tree of a free action on
an arbitrary nonempty set, as well as its extension to
a simplicial tree on a group naturally associated to
the given free action, in order to use it further for
obtaining by suitable deformations more general arboreal
structures. 

Let $H$ be a group acting freely on a nonempty set $X$. Let 
$B = \{b_i\,|\,i \in I\} \sse X$ be a set of representatives
for the $H$-orbits. The bijection $H \times I \lra X, (h, i) 
\mapsto h b_i$ identifies up to isomorphism the $H$-set $X$
with the cartesian product $H \times I$ with the canonical
free action of the group $H$, $H \times (H \times I) \lra H \times I,
(h_1, (h_2, i)) \mapsto (h_1 h_2, i)$. 

We assume that $I \cap H = \{1\}$, and we shall take
$b_1 = (1, 1)$ as {\em basepoint} in $X \cong H \times I$. 
Set $I' := I \sm \{1\}, X' := \sqcup_{i \in I'} Hb_i \cong H \times I'$. 

\subsection{The underlying tree of the $H$-set $X$}

$\quad$ The set $X \cong H \times I$ has a natural structure of
simplicial tree with the elements of $X$ as vertices, and 
the ordered pairs $(b_1, h b_1), h \in H \sm \{1\}$, and
$(h b_1, h b_i), h \in H, i \in I'$, as oriented edges.
Taking $b_1$ as root, we obtain a rooted order-tree
$(X, b_1, \leq)$ with the partial order $\leq$ given by
$$x < y \Llra\, {\rm either}\, x = b_1, y \neq b_1\,
{\rm or}\, \exists h \in H \sm \{1\}, i \in I', x = h b_1, y = h b_i,$$
and the induced meet-semilattice operation $\wedge$ and 
locally linear median operation 
$$Y(x, y, z) = (x \wedge y) \vee (y \wedge z) \vee
(z \wedge x) \in \{x \wedge y, y \wedge z, z \wedge x\}.$$

\begin{rems} \em
$(1)$ With respect to the partial order $\leq$, $b_1$ is
the least element of $X$, while $X'$ is the set of all maximal
elements of $X$. In particular, for $x, y \in X, 
x \wedge y \in X' \Llra x = y \in X'$, and hence 
$x \wedge y \in H b_1$ provided $x \neq y$, whence
$Y(x, y, z) \in H b_1$ whenever $x \neq y, y \neq z, z \neq x$.

$(2)$ $H b_1$ is a retractible convex subset of the median
set $(X, Y)$, with the associated folding 
$\theta : X \lra X, h b_i \mapsto h b_1$, compatible with 
the action of $H$, i.e. $\theta(h x) = h \theta(x)$ for
all $h \in H, x \in X$.

$(3)$ The median operation $Y$ is {\em almost compatible}
with the action of $H$ on $X$ in the following sense: 
$Y(h x, h y, h z) \in H\, Y(x, y, z)$ for all $h \in H, 
x, y, z \in X$, i.e. the map $Y : X^3 \lra X$ induces 
a map $H \sm X^3 \lra H \sm X$. Indeed, for 
$h \in H, x, y, z \in X$, we obtain
\begin{center}
$Y(h x, h y, h z) = \left\{ \begin{array}{lcl}
h Y(x, y, z)  & \mbox{if} & |\{\theta(x), \theta(y), \theta(z)\}| \in \{1, 2\}  \\
Y(x, y, z) = b_1 & \mbox{if} & |\{\theta(x), \theta(y), \theta(z)\}| = 3
\end{array}
\right.$
\end{center}
Consequently, the necessary and sufficient condition for
the median operation $Y$ to be compatible
with the action of $H$ on $X$, i.e.
$h Y(x, y, z) = Y(h x, h y, h z)$ for all $h \in H, x, y, z \in X$,
is that either $H = 1$ or $H \cong \Z/2$.
\end{rems}

\subsection{The group $\wh{H}$ and its underlying tree}

$\quad$ We denote by $F$ the free group with base $I'$, and by 
$\wh{H} := H \ast F$ the free product of the  
groups $H$ and $F$. The group $H$ is canonically identified with 
a subgroup of $\wh{H}$, while the injective map $\iota : X \lra
\wh{H}, h b_i \mapsto h i$, identifies the $H$-set 
$X \cong H \times I$ with the disjoint union 
$H \bigsqcup (\bigsqcup_{i \in I'} H i) \sse \wh{H}$ on
which $H$ acts freely by left multiplication. 

Using the natural tree structure of the free group $F$, 
we extend as follows the underlying tree of the $H$-set $X$ 
as defined in 2.1. to a simplicial tree on the 
underlying set of the free product $\wh{H} = H \ast F$.

Let $l : \wh{H} \lra \N$ denote the length function associated
to the system of generators $J = J^{- 1} := (H \sm \{1\}) \sqcup I'^{\pm 1}$, 
so $l(w)$ is the minimum length of any expression $w = w_1 \cdots w_n$
with $w_k \in J$. In particular, $l(w) = 0 \Llra w = 1$, $l(w^{- 1}) = l(w)$ 
for all $w \in \wh{H}$, and $l(u v) \leq l(u) + l(v)$ for all 
$u, v \in \wh{H}$. Since $\wh{H} = H \ast F$
and $F$ is free with base $I'$, it follows that the
expression of minimal length above is unique for any $w \in \wh{H}$;
call it the {\em reduced normal form} of $w$, and set o$(w) := w_1$,
t$(w) := w_n$ provided $l(w) = n \geq 1$.
For $u, v \in \wh{H}$, put $u \leq v \Llra
l(v) = l(u) + l(u^{- 1} v)$, and write $v = u \bullet (u^{- 1} v)$
provided $u \leq v$. The binary relation $\leq$ is a partial order
extending the partial order $\leq$ on $X$ as defined in 2.1. Moreover
the partial order $\leq$ makes $\wh{H}$ a {\em rooted order-tree} 
with $1$ as distinguished base point, the {\em root}, while the 
corresponding meet-semilattice operation $\wedge$ and locally 
linear median operation $Y$ are extensions
of the operations $\wedge$ and $Y$ on $X$. 

For all $u, v \in \wh{H}$, the cell 
$[u, v] := \{Y(u, v, w)\,|\,w \in \wh{H}\}$ is the union
of the closed intervals $[u \wedge v, u]$ and 
$[u \wedge v, v]$. Since the cell $[u, v]$ has finitely many
elements for all $u, v \in \wh{H}$, $\wh{H}$ is a
$\Z$-tree with the distance function
$d : \wh{H} \times \wh{H} \lra \Z$ defined by 
$d(u, v) := |[u, v]| - 1 = l(w^{- 1} u) + l(w^{- 1} v)$, where $w = u \wedge v$.
Thus $d(u, v) = l(u^{- 1} v)$ provided $u \leq v$, in particular,
$d(u, 1) = l(u)$ for all $u \in \wh{H}$, 
and hence $d(u, v) = d(u, u \wedge v) + d(v, u \wedge v) \geq
l(u^{- 1} v)$ for all $u, v \in \wh{H}$. Consequently,
$d(u, v) = l(u^{- 1} v) \Llra u^{- 1} v = (u^{- 1} (u \wedge v))
\bullet ((u \wedge v)^{- 1} v)$, while $d(u, v) = l(u^{- 1} v) +1$
otherwise. The latter situation holds whenever $u \neq u \wedge v
\neq v$ and o$((u \wedge v)^{- 1} u)$, 
o$((u \wedge v)^{- 1} v) \in H \sm \{1\}$.
In graph theoretic terms, the underlying tree of $\wh{H}$
has the elements of $\wh{H}$ as vertices, and the ordered
pairs $(u, v)$, with $u \leq v, l(u^{- 1} v) = 1$, as oriented
edges.
\smallskip

$X$ is a retractible convex subset of 
the median set $(\wh{H}, Y)$, with the canonical retract 
$\varphi : \wh{H} \lra X$ defined by 
$\varphi(w) :=$ the greatest element $x \in X$ for which 
$x \leq w$, i.e. $w = x \bullet (x^{- 1} w)$. Thus $\varphi(w)
= h \in H \Llra w = h$ or $h i^{- 1} \leq w$ for some $i \in I'$, 
while $\varphi(w) = h i$ with $h \in H, i \in I' \Llra h i \leq w$.
In particular, $H = H b_1$ is a retractible convex subset 
of the median set $(\wh{H}, Y)$ with the retract 
$\wh{\theta} := \theta \circ \varphi : \wh{H} \lra H$. 
Thus $\wh{\theta}(w) = 1 \Llra\,$ either $w = 1$ or o$(w) \in I'^{\pm 1}$,
and $\wh{\theta}(w) = h \in H \sm \{1\} \Llra\,$ o$(w) = h$.

\begin{rem} \em
$(1)$ The maps $\varphi$ and $\wh{\theta}$ are $H$-equivariant.

$(2)$ Let $u, v \in \wh{H}$ be such that $u^{- 1} \wedge v = 1$.
Then the following assertions are equivalent.

(i) $\varphi(u v) = \varphi(u)$.

(ii) Either $u \not \in H$ or $\varphi(v) = 1$.

$(3)$ $\varphi(u x) = \varphi(u)$ for all 
$u \in \wh{H} \sm H, x \in X$.
\end{rem}

The next lemma collects some basic properties which describe the
relation between the group $\wh{H}$ and its underlying
tree.

\begin{lem} 
\ben
\item[\rm (1)] $u \leq v$ and $v \leq w \Lra w^{- 1} v \leq w^{- 1} u$.

\item[\rm (2)] For $u, v \in \wh{H}$, let $a := u^{- 1} \wedge v,
u' := u a, v' := a^{- 1} v, u'' := u' \wh{\theta}(u'^{- 1}),
v'' := \wh{\theta}(v')^{- 1} v', h := \wh{\theta}(u'^{- 1})^{- 1} 
\wh{\theta}(v') \in H$, with $h = 1 \Llra \wh{\theta}(u'^{- 1}) =
\wh{\theta}(v') = 1$. Then $u v = u'' h v'' = u'' \bullet h \bullet v''$,
i.e. $l(u v) = l(u'') + l(h) + l(v'')$. 

In particular, $u \leq v \Llra u^{- 1} \wedge (u^{-1} v) = 1$ and
either $\wh{\theta}(u^{- 1}) = 1$ or $\wh{\theta}(u^{- 1} v) = 1$.

\item[\rm (3)] The necessary and sufficient condition for $\wh{H}$
together with the median operation $Y$ to be a median group is that
either $H = 1$ or $H \cong \Z/2$.

\item[\rm (4)] For all $s, u, v \in \wh{H}, s u \wedge s v \leq
s Y(u, v, s^{- 1})$, with $(s u \wedge sv)^{- 1} s Y(u, v, s^{- 1}) \in H$.
 
\item[\rm (5)] For all $u, v, w \in \wh{H}, Y(t^{- 1} u, t^{- 1} v, t^{- 1} w) = 1$, 
where $t := Y(u, v, w)$.
\een
\end{lem}

\bp
The proof of the assertions (1) and (2) is straightforward.

(3). It follows by Proposition 1.11. that the necessary and
sufficient condition for $(\wh{H}, Y)$ to be a median group
is that $u^{- 1} (u \wedge v) \leq u^{- 1} v$ for all $u, v \in \wh{H}$.
According to (2), the latter condition holds if and only if
for all $u, v \in \wh{H}$, either $\wh{\theta}((u \wedge v)^{- 1} u) = 1$
or $\wh{\theta}((u \wedge v)^{- 1} v) = 1$. One checks easily
that the last sentence is equivalent with $|H| \leq 2$.
\smallskip

(4). Let $s, u, v \in \wh{H}$. Setting $a := u \wedge v \wedge s^{- 1},
b := Y(u, v, s^{- 1})$, we have to show that 
$s b = (s u \wedge sv) \bullet p$ with $p \in H$. We distinguish 
the following three cases:
\smallskip

(4.1.) $a = u \wedge s^{- 1} = v \wedge s^{- 1} \leq b = u \wedge v$: 
We may assume without loss that $a = 1$ since
$s = s' \bullet a^{- 1}, u = a \bullet u', v = a \bullet v'$ 
with $a' := u' \wedge s'^{- 1} = v' \wedge s'^{- 1} = 1$, 
therefore $s u = s' u', sv = s' v'$, and 
$$s Y(u, v, s^{- 1}) = s (u \wedge v) = s' (u' \wedge v') =
s' Y(u', v', s'^{- 1}),$$
so we may replace the elements $u, v, s$ by $u', v', s'$ respectively.
As $u \wedge s^{- 1} = 1$, it follows that either 
$s u = s \bullet u$ or $s u = (s' \bullet g)
(h \bullet u') = s' \bullet (g h) \bullet u'$ where 
$g, h \in H \sm \{1\}, g h \neq 1$. A similar alternative 
holds for the pair $(s, v)$. Thus we have the following possible 
situations:
\smallskip

(4.1.1.) $s u = s \bullet u, s v = s \bullet v$: Then
$s u \wedge s v = s \bullet b$, so $p = 1$.

(4.1.2.) $s u \neq s \bullet u, s v = s \bullet v$: Then
$s = s' \bullet g, u = h \bullet u'$ with $g, h \in H \sm \{1\},
g h \neq 1$, therefore $s u = s' \bullet (g h) \bullet u'$,
and either $v = 1$ or o$(v) \in I'^{\pm 1}$. Consequently, 
$b = 1, s u \wedge sv = s', s b = s = s' \bullet p$ with
$p = g \in H \sm \{1\}$.

(4.1.3.) $s u = s \bullet u, s v \neq s \bullet v$: We proceed
as in case (4.1.2).

(4.1.4.) $s u \neq s \bullet u, s v \neq s \bullet v$: Then
$s = s' \bullet g, u = h_1 \bullet u', v = h_2 \bullet v'$
with $g, h_j \in H \sm \{1\}, g h_j \neq 1, j = 1, 2$. We have
the alternative: either $h := h_1 = h_2 $ or $h_1 \neq h_2$. In
the first case we obtain  $b = h \bullet (u' \wedge v'),
sb = s u \wedge s v = s' \bullet (g h) \bullet (u' \wedge v')$,
so $p = 1$, while in the second case we get $b = 1, s u \wedge
s v = s', s b = s = s' \bullet p$ with $p = g \in H \sm \{1\}$
as desired. 
\smallskip

(4.2.) $a = u \wedge s^{- 1} = u \wedge v < b = v \wedge s^{- 1}$: 
As in case (4.1.), we may assume that $a = 1$. We get
$s = s' \bullet b^{- 1}, v = b \bullet v'$ with 
$u \wedge b = v' \wedge s'^{- 1} = 1$. We have to show that
$s' = (s' b^{- 1} u \wedge s' v') \bullet p$ with $p \in H$.

We have the following possible situations:
\smallskip

(4.2.1.) $b^{- 1} u = b^{- 1} \bullet u, s' v' = s' \bullet v'$:
Then $s' b^{- 1} u \wedge s' v' = s'$, so $p = 1$.

(4.2.2.) $b^{- 1} u = b^{- 1} \bullet u, s' v' \neq s' \bullet v'$:
Then $s' = s'' \bullet g, v' = h \bullet v''$ with 
$g, h \in H \sm \{1\}, g h \neq 1$. The desired
result follows with $p = g \in H \sm \{1\}$.

(4.2.3.) $b^{- 1} u \neq b^{- 1} \bullet u, s' v' = s' \bullet v'$:
Then $b = g \bullet b', u = h \bullet u'$ with $g, h \in H \sm \{1\},
g \neq h$. First let us assume that $b' \neq 1$. Since 
$s = s' \bullet b^{- 1}$, it follows that
$s' b^{- 1} u = s' \bullet b'^{- 1} \bullet (g^{- 1} h) \bullet u'$. As
$v = b \bullet v' = g \bullet b' \bullet v'$, we deduce that 
$v' \wedge b'^{- 1} = 1$. Since $s' v' = s' \bullet v'$ by 
assumption, the required result follows with $p = 1$.

Next let us assume that $b' = 1$, i.e. $b = g \in H \sm \{1\}$.
Then $b^{- 1} u = (g^{- 1} h) \bullet u', s = s' \bullet b^{- 1} =
s' \bullet g^{- 1}$, therefore either $s' = 1$ or t$(s') \in I'^{\pm 1}$.
Consequently, $s' b^{- 1} u = s' \bullet (g^{- 1} h) \bullet u'$.
On the other hand, since $v = b \bullet v' = g^{- 1} \bullet v'$,
it follows that either $v' = 1$ or o$(v') \in I'^{\pm 1}$, and
hence $v' \wedge (g^{- 1} h) = 1$. As $s' v' = s' \bullet v'$,
we get as above the desired result with $p = 1$.

(4.2.4.) $b^{- 1} u \neq b^{- 1} \bullet u, s' v' \neq s' \bullet v'$:
According to (4.2.3.) we get $s' \leq s' b^{- 1} u$. On the other hand,
it follows by assumption that $s' = s'' \bullet g', v' = h' \bullet v''$ 
with $g', h' \in H \sm \{1\}, g' h' \neq 1$. Thus $s' v' \wedge s' =
(s'' \bullet (g' h') \bullet v'') \wedge (s'' \bullet g') = s'' < s'$,
and hence the required result with $p = g' \in H \sm \{1\}$.
\smallskip

(4.3) $a = v \wedge s^{- 1} = u \wedge v < b = u \wedge s^{- 1}$:
We proceed as in case (4.2.).
\smallskip

(5). Let $u, v, w \in \wh{H}, t := Y(u, v, w)$. Since $Y(u, v, t) = t$,
it follows by (4) that $t^{- 1} u \wedge t^{- 1} v \leq t^{- 1} Y(u, v, t) = 1$,
therefore $t^{- 1} u \wedge t^{- 1} v = 1$, and similarly,
$t^{- 1} v \wedge t^{- 1} w = t^{- 1} w \wedge t^{- 1} u = 1$,
and hence $Y(t^{- 1} u, t^{- 1} v, t^{- 1} w) = 1$ as desired. 
\ep


\section{The deformation of the underlying tree of $\wh{H}$ into
median group operations}

$\quad$ Let $H$ be a group acting freely on a nonempty set $X$.
As shown in Section 2, the free action $H \times X \lra X$
is extended through the $H$-equivariant embedding $\iota : X \lra \wh{H}$, 
with the $H$-equivariant retract $\varphi : \wh{H} \lra X$,
to the transitive and free action of the group $\wh{H} = H \ast F$ 
on itself given by left multiplication. Thus 
the conditions (1) and (2) of Theorem 1 (see Introduction) 
are obviously satisfied.

We denote by $\cM(X)$ the set of all median operations $m$ on
$X$ which are compatible with the action of $H$, while
by $\cM(\wh{H}, \varphi)$ we denote the set of those 
median group operations $\wh{m}$ on $\wh{H}$ for which
the retract $\varphi$ is a folding identifying $X$ with 
a retractible convex subset of the median set $(\wh{H}, \wh{m})$. 
We denote by $\cM_l(X), \cM_l(\wh{H}, \varphi)$ ($\cM_s(X),
\cM_s(\wh{H}, \varphi)$) the subsets of $\cM(X)$ and
$\cM(\wh{H}, \varphi)$ respectively consisting of those median
operations which are locally linear (simplicial).
According to Lemma 1.22., the restriction
map ${\rm res} : \cM(\wh{H}, \varphi) \lra \cM(X)$ 
is injective, whence the induced maps 
${\rm res}_l : \cM_l(\wh{H}, \varphi) \lra \cM_l(X)$
and ${\rm res}_s : \cM_s(\wh{H}, \varphi) \lra \cM_s(X)$
are injective too.

The present section is devoted to the proof of a more 
explicit version of Theorem 1. With the notation above we obtain

\begin{te} The map ${\rm res} : \cM(\wh{H}, \varphi) \lra \cM(X)$ 
is bijective. Let $m \in \cM(X)$. Then the following assertions hold.
\ben
\item[\rm (1)] The unique median group operation
$\wh{m} \in \cM(\wh{H}, \varphi)$ lying over $m$ is a deformation 
of the underlying simplicial tree of $\wh{H}$ induced by the median 
operation $m$ and the retract $\varphi$, defined by
$$\wh{m}(u, v, w) = t\, m(\varphi(t^{- 1} u), \varphi(t^{- 1} v), 
\varphi(t^{- 1} w))$$
for $u, v, w \in \wh{H}$, where $t = Y(u, v, w)$.
\item[\rm (2)] The induced meet-semilattice operation
$u \cap v := \wh{m}(u, 1, v)$ is a deformation of the meet-semilattice
operation $\wedge$, defined by \em
$$u \cap v = (u \wedge v) m(\varphi((u \wedge v)^{- 1} u), \varphi((u \wedge v)^{- 1}),
\varphi((u \wedge v)^{- 1} v))$$
for $u, v \in \wh{H}$. \em
\item[\rm (3] The induced partial order $\subset$ is defined by \em
$$u \subset v \Llra u \cap v = u \Llra (u \wedge v)^{- 1} u \in \texttt{[}
\varphi((u \wedge v)^{- 1}), \varphi((u \wedge v)^{- 1} v)\,
\texttt{]} \sse X.$$ \em
In particular, $u \subset v$ provided $u \leq v$ and 
{\em $1 \in \texttt{[} \varphi(u^{- 1}), \varphi(u^{- 1} v) \texttt{]}$}. 
\een
\end{te}

\bp \footnote{$^)$ See \cite[3, Theorem 3.1.]{v2} for an alternative 
more technical proof}
Let $m \in \cM(X)$. We have to define a median group operation
$\wh{m} \in \cM(\wh{H}, \varphi)$ such that $\wh{m}(x, y, z) =
m(x, y, z)$ for all $x, y, z \in X$. Since $m$ is compatible
with the action of $H$, the group $H$ acts from the right on
$Spec(X, m)$ according to the rule 
$\fp^h := h^{- 1} \fp = \{h^{- 1} x\,|\,x \in \fp\}$ for 
$\fp \in Spec(X, m), h \in H$. 

On the other hand, we consider 
the natural action from the right of the group $\wh{H}$ on the 
power set $\cP(\wh{H})$, $(P, u) \mapsto P^u := u^{- 1} P = 
\{u^{- 1} v\,|\,v \in P\}$. Let 
$$\cS := \{\varphi^{- 1}(\fp)^u =
u^{- 1} \varphi^{- 1}(\fp)\,|\,\fp \in Spec(X, m), u \in \wh{H}\}$$
be the $\wh{H}$-orbit of the subset $\varphi^{- 1}(Spec(X, m))$. 

The set $\cS$ is closed under the involution
$P \mapsto \wh{H} \sm P$ which is compatible with the action of $\wh{H}$,
i.e. $\wh{H} \sm (P^u) = (\wh{H} \sm P)^u$ for all $P \in \cS, u \in \wh{H}$. 

Let $P = \varphi^{- 1}(\fp)^u \in \cS \sm 
\{\es, \wh{H}\}$, whence $\fp \neq \es, X$. Then $P \cap X = 
\{x \in X\,|\,\varphi(u x) \in \fp\}$.
We distinguish the following three cases.

(i) $u \in H$. Then $P \cap X = \fp^u \in Spec(X, m) \sm \{\es, X\}$
and $P = \varphi^{- 1}(\fp^u) \in \varphi^{- 1}(Spec(X, m) \sm \{\es, X\})$.

(ii) $u \not \in H$ and $\varphi(u) \in \fp$. Then $\varphi(u x) = \varphi(u) 
\in \fp$ for all $x \in X$, by Remark 2.2.(3), whence $X \sse P$.

(iii) $u \not \in H$ and $\varphi(u) \in X \sm \fp$. Then $X \sse \wh{H} \sm P$
by (ii), and hence $P \cap X = \es$.

Consequently, the $H$-equivariant retract $\varphi : \wh{H} \lra X$ of
the $H$-equivariant embedding $X \lra \wh{H}$ induces a $H$-equivariant
embedding $Spec(X, m) \lra \cS, \fp \mapsto \varphi^{- 1}(\fp)$ with
the $H$-equivariant retract $\cS \lra Spec(X, m), P \mapsto P \cap X$,
and $\varphi^{- 1}(Spec(X, m)) = \{\es, \wh{H}\} \cup \{P \in \cS\,|\,
P \cap X \not \in \{\es, X\}\}$.

$\wh{H}$ acts canonically on the power set $\cP(\cS)$ according to the
rule 
$$u\, \cU := \{P^{u^{- 1}} = u P\,|\,P \in \cU\}\,{\rm for}\, u \in \wh{H},
\cU \sse \cS.$$
Define the negation operator on $\cP(\cS), \cU \mapsto
\neg \cU := \{P \in \cS\,|\,\wh{H} \sm P \not \in \cU\}$. It follows that
$\neg \cU \sse \neg \cV \Llra \cV \sse \cU, \neg (\cU \cup \cV) =
(\neg \cU) \cap (\neg \cV), \neg (\cU \cap \cV) = (\neg \cU) \cup (\neg \cV)$,
and $\neg (\neg \cU) = \cU$ for $\cU, \cV \sse \cS$. In addition,
the operator $\neg$ is compatible with the action of $\wh{H}$, i.e.
$\neg (u\, \cU) = u (\neg \cU)$ for all $u \in \wh{H}, \cU \sse \cS$.

The group $\wh{H}$ and the power set $\cP(\cS)$ are also related 
through the map 
$$\wh{H} \lra \cP(\cS), u \mapsto U(u) := \{P \in \cS\,|\,
u \not \in P\}.$$ 
Notice that $\neg U(u) = U(u)$ and $u U(v) = U(u v)$ for
$u, v \in \wh{H}$, so $\wh{H}$ acts transitively on the $\neg$-invariant set 
$\{U(u)\,|\,u \in \wh{H}\}$. To show that the action is free, let
$u \in \wh{H} \sm \{1\}$.  It suffices to show that $U(u) \not \sse U(1)$.
We distinguish the following two cases.

(i) $\varphi(u) \neq 1$. Then, by Theorem 1.6.(1), there is 
$\fp \in Spec(X, m)$ such that $1 \in \fp, \varphi(u) \not \in \fp$,
and hence $\varphi^{- 1}(\fp) \in U(u) \sm U(1)$.

(ii) $\varphi(u) = 1$. Then, since $u \neq 1$, it follows that
$u = i^{- 1} \bullet v$ for some $i \in I', v \in \wh{H}$ with
$\varphi(v) \neq i$. Consequently, by Theorem 1.6.(1) again,
there is $\fp \in Spec(X, m)$ such that $i \in \fp, \varphi(v) \not
\in \fp$, therefore $\varphi^{- 1}(\fp)^i \in U(u) \sm U(1)$ as required.

To obtain the desired median group operation $\wh{m}$ on $\wh{H}$,
identified through the injective map $u \mapsto U(u)$ with a 
$\wh{H}$-subset of $\cP(\cS)$, we have to show that the
subset $\{U(u)\,|\,u \in \wh{H}\}$ is closed under the canonical
median operation $\fm$ on $\cP(\cS)$ 
$$\fm(\cU, \cV, \cW) := (\cU \cap \cV) \cup (\cV \cap \cW) \cup
(\cW \cap \cU) = (\cU \cup \cV) \cap (\cV \cup \cW) \cap (\cW \cup \cU),$$
which is compatible with the action of $\wh{H}$, i.e. $u\, \fm(\cU, \cV, \cW) =
\fm(u\, \cU, u \cV, u \cW)$ for $u \in \wh{H}, \cU, \cV, \cW \sse \cS$.

Taking into acount the invariance of the $U(u)$'s under the negation 
operator $\neg$, it suffices to show that for arbitrary $u, v, w \in \wh{H}$,
$\fm(U(u), U(v), U(w)) \sse U(t x)$, where $t := Y(u, v, w), x := 
m(\varphi(t^{- 1} u), \varphi(t^{- 1} v), \varphi(t^{- 1} w))$.
Using the
compatibility of the median operation $\fm$ with 
the action of $\wh{H}$ and the fact that $Y(t^{- 1} u, t^{- 1} v, t^{- 1}w) = 1$
by Lemma 2.3.(5), it remains to show that $\fm(U(u), U(v), U(w)) \sse U(x)$,
where $x := m(\varphi(u), \varphi(v), \varphi(w))$, for $u, v, w \in \wh{H}$
satisfying $Y(u, v, w) = 1$. Consequently, it suffices to show that
$U(u) \cap U(v) \sse U(x)$ provided $u \wedge v = 1$ and $x$ belongs to
the cell $\texttt{[} \varphi(u), \varphi(v) \texttt{]}$ of the median set $(X, m)$
\footnote{$^)$ For $x, y \in X$, the cell $\texttt{[} x, y \texttt{]}$ of 
the median set $(X,m)$ must not be confounded with the cell $[x, y]$ 
of the median set $(X, Y)$ as defined in 2.1.}.
Assuming the contrary, let $P \in \cS$ be such that $x \in P, u, v \not \in P$,
in particular, $P \neq \wh{H}$ and $P \cap X \neq \es$. 
We distinguish the following two cases.
\smallskip

(i) $\fp := P \cap X \neq X$, whence $P = \varphi^{- 1}(\fp)$.
As $x \in \texttt{[} \varphi(u), \varphi(v) \texttt{]} \cap \fp$, it
follows that either $\varphi(u) \in \fp$ or $\varphi(v) \in \fp$, whence
either $u \in P$ or $v \in P$, and hence a contradiction.

(ii) $X \sse P$, whence $P = \varphi^{- 1}(\fq)^s$ for some
$\fq \in Spec(X, m) \sm \{\es, X\}, s \in \wh{H} \sm H$ with 
$\varphi(s) \in \fq$. Since $u \wedge v = 1$ by assumption,
it follows that either $s^{- 1} \wedge u = 1$ or $s^{- 1} \wedge v = 1$,
therefore either $\varphi(s u) = \varphi(s) \in \fq$ or
$\varphi(s v) = \varphi(s) \in \fq$ according to Remark 2.2.(2).
Consequently, either $u \in P$ or $v \in P$, again a contradiction.

Thus we have obtained the desired median group operation 
$\wh{m} : \wh{H}^3 \lra \wh{H}$, inducing the meet-semilattice operation
$\cap$ and the partial order $\subset$ as defined in the statements (2) 
and (3) of the theorem.
One checks easily that res$(\wh{m}) = m$ and $\varphi$ is a
folding of the median set $(\wh{H}, \wh{m})$ as required.
According to Lemma 1.22., $\wh{m}$ is unique with the properties 
above, and $Spec(\wh{H}, \wh{m}) = \cS$. This completes the proof.
\ep

Tha next lemma provides equivalent descriptions for the partial
order $\subset$ \footnote{$^)$ For other descriptions of the 
partial order $\subset$ see \cite[Lemma 3.5.]{v2}}.

\begin{lem}
Let $m \in \cM(X)$ and $\subset$ be the partial order induced
by the unique median group operation $\wh{m} \in \cM(\wh{H}, \varphi)$
lying over $m$. Then the following assertions are equivalent
for $u, v \in \wh{H}$.
\ben
\item[\rm (1)] $u \subset v$.
\item[\rm (2)] $U(1) \cap U(v) \sse U(u)$, where $U(w) := \{P \in \cS\,|\,
w \not \in P\}$ for $w \in \wh{H},\, \cS := Spec(\wh{H}, \wh{m})$.
\item[\rm (3)] There is $w \in \wh{H}$ such that $w \leq v$
and {\em $w^{- 1} u \in \texttt{[} \varphi(w^{- 1}), \varphi(w^{- 1} v)
\texttt{]} \sse X$}. 
\item[\rm (4)] Either $u = v = 1$ or there is $w \in \wh{H}$ such that \em
$$w \leq v, \varphi(w^{- 1}) \in I, \varphi(w) = 1 \Lra \varphi(w^{- 1} v) \neq 1,
w^{- 1} u \in \texttt{[} \varphi(w^{- 1}), \varphi(w^{- 1} v) 
\texttt{]} \sse X,$$ 
\em whence $w \varphi(w^{- 1}) \subset u \subset w \varphi(w^{- 1} v)
\subset v$, while $w \varphi(w^{- 1}) < w \bullet \varphi(w^{- 1} v) \leq v$.   
\een 
\end{lem}

\bp
(1) $\Llra$ (2) follows by Theorem 1.6.(1), while (4) $\Lra$ (3) is obvious.

(1) $\Lra$ (4). Let $u, v \in \wh{H}$ be such that $u \subset v$.
If $v = 1$ then $u = 1$, so let us assume that $v \neq 1$. Setting
$a := u \wedge v, u = a \bullet b, v = a \bullet c$ with $b \wedge c = 1$,
we have by assumption $b \in \texttt{[} \varphi(a^{- 1}), \varphi(c) \texttt{]}$. 
We distinguish the following three cases.

(i) $\varphi(a^{- 1}) \in I$, with $\varphi(c) \neq 1$ provided 
$\varphi(a^{- 1}) = 1$. Then $w := a$ satisfies the requirements.

(ii) $\varphi(a^{- 1}) = \varphi(c) = 1$. Then $b = 1, u = a \leq v = u \bullet c$.
As $v \neq 1$ it follows that either $u \neq 1$ or $c \neq 1$.

First assume that $u \neq 1$, whence $u = u' \bullet i$ with $i \in I'$ 
since $\varphi(u^{- 1}) = 1, u \neq 1$ by assumption. Then $w := u'$ 
satisfies the required conditions provided $\varphi(u'^{- 1}) \in I$.
Assuming that $\varphi(u'^{- 1}) \not \in I$, we obtain $u = u'' \bullet
h \bullet i$ with $h \in H \sm \{1\}, \varphi(u''^{- 1}) \in I$, therefore
$w := u''$ satisfies the requirements.

Next assume that $c \neq 1$, whence $c = i^{- 1} \bullet c'$ with $i \in I'$
since $\varphi(c) = 1, c \neq 1$. Then $w := u \bullet i^{- 1}$ satisfies
the required conditions.

(iii) $\varphi(a^{- 1}) \not \in I$, whence $a = w \bullet h$ with
$h \in H \sm \{1\}, \varphi(w^{- 1}) \in I$. Then
$\varphi(w^{- 1} v)  = \varphi(h \bullet c) = h \bullet \varphi(c) \neq 1$,
and $w^{- 1} u = h \bullet b \in h \texttt{[} \varphi(a^{- 1}), \varphi(c) 
\texttt{]} = \texttt{[} \varphi(w^{- 1}), \varphi(w^{- 1} v) \texttt{]}$ as
desired. 

One checks easily that in all the cases above, 
$w \varphi(w^{- 1}) < w \varphi(w^{- 1} v) \leq v$ and
$w \varphi(w^{- 1}) \subset u \subset w \varphi(w^{- 1} v) \subset v$.

(3) $\Lra$ (2). Let $u, v \in \wh{H}$ be such that $w^{- 1} u \in
\texttt{[} \varphi(w^{- 1}), \varphi(w^{- 1} v) \texttt{]} \sse X$ for some
$w \leq v$. Assuming that $U(1) \cap U(v) \not \sse U(u)$, let 
$P \in \cS$ be such that $1, v \not \in P, u \in P$, whence
$w^{- 1}, w^{- 1} v \not \in Q, w^{- 1} u \in Q$, where $Q := P^w$. 
As $w^{- 1} u \in X$ by assumption, it follows that $Q \cap X \neq \es$,
and hence there are the following two possibilities.

(i) $\fq := Q \cap X \neq X$, whence $Q = \varphi^{- 1}(\fq)$. 
As $w^{- 1}, w^{- 1} v \not \in Q$, we obtain $\varphi(w^{- 1}),
\varphi(w^{- 1} v) \in X \sm \fq$, therefore, by Theorem 1.6.(1), 
$\fq\, \cap \texttt{[} \varphi(w^{- 1}), \varphi(w^{- 1} v) \texttt{]} = \es$,
contrary to the assumption that $w^{- 1} u \in \fq \cap \texttt{[}
\varphi(w^{- 1}, \varphi(w^{- 1} v) \texttt{]}$.

(ii) $X \sse Q$, whence $Q = \varphi^{- 1}(\fq)^s$ for some 
$\fq \in Spec(X, m) \sm \{\es, X\}, s \in \wh{H} \sm H$ with 
$\varphi(s) \in \fq$. As $w^{- 1}, w^{- 1} v \not \in Q$, we obtain
$\varphi(s w^{- 1}), \varphi(s w^{- 1} v) \in X \sm \fq$.
On the other hand, since $w \leq v$, it follows that
$w^{- 1} \wedge w^{- 1} v = 1$, and hence either 
$s^{- 1} \wedge w^{- 1} = 1$ or $s^{- 1} \wedge w^{- 1} v = 1$.
Consequently, by Remark 2.2.(2), we deduce that either
$\varphi(s w^{- 1}) = \varphi(s) \in \fq$ or $\varphi(s w^{- 1} v) =
\varphi(s) \in \fq$, i.e. a contradiction. This finishes the proof.
\ep

To any $v \in \wh{H} \sm \{1\}$ we associate the following two sets
$$C_v := \{w \in \wh{H}\,|\,w \leq v, \varphi(w^{- 1}) \in 
I, \varphi(w^{- 1}) = 1 \Lra \varphi(w^{- 1} v) \neq 1\}$$ 
and
$$O_v := \{w \varphi(w^{- 1})\,|\,w \in C\},$$
together with the map $\zeta : C_v \lra O_v$ defined by 
$\zeta(w) = w \varphi(w^{- 1})$, whence $\zeta(w) \leq w$ for $w \in C_v$. 
$C_v$ and $O_v$ are nonempty finite sets, totally ordered with respect 
to $\leq$, and the map $\zeta$ is an isomorphism of totally ordered sets.
Setting $C_v = \{w_i\,|\,i = \overline{1,\,n}\,\}$ with $n \geq 1, w_i <
w_{i + 1}$, it follows that $\zeta(w_1) = 1, w_i \leq \zeta(w_{i + 1}) = w_i \bullet
\varphi(w_i^{- 1} v)$ for $i < n$, and $\zeta(w_n) < v = w_n \bullet \varphi(w_n^{- 1} v)$.
Thus the totally ordered finite set $([1, v], \leq)$ is the union of
$n$ adjacent proper closed intervals $J_i := [\zeta(w_i), \zeta(w_{l + 1}],
i = \overline{1,\,n - 1}, J_n := [\zeta(w_n), v]$ with $w_i \in J_i, i = \overline{1,\,n}$.
Call this decomposition in adjacent closed intervals the {\em combinatorial
configuration associated to the element} $v \in \wh{H} \sm \{1\}$.

For instance, taking $v = i^{- 2} h j^3 g$ with $h, g \in H \sm \{1\},
i, j \in I', l(v) = 7$, we obtain the totally ordered sets
$$C_v = \{i^{- 1}, i^{- 2}, i^{- 2} h j, i^{- 2} h j^2, i^{- 2} h j^3\},\,
O_v = \{1, i^{- 1}, i^{- 2} h j, i^{- 2} h j^2, i^{- 2} h j^3\}$$
of cardinality $n = 5$, and the adjacent closed intervals
$$J_1 = [1, i^{- 1}], I_2 = [i^{- 1}, i^{- 2} h j], J_3 = [i^{- 2} h j,
i^{- 2} h j^2], J_4 = [i^{- 2} h j^2, i^{- 2} h j^3], J_5 = [i^{- 2} h j^3, v]$$
of cardinality $2, 4, 2, 2, 2$ respectively.

As a consequence of Lemma 3.2. and Corollaries 1.19., 1.20., we obtain

\begin{co}
Let $m \in \cM(X)$ and $\subset$ be the partial order induced
by the unique median group operation $\wh{m} \in \cM(\wh{H}, \varphi)$
lying over $m$. Then, for any $v \in \wh{H} \sm \{1\}$, the cell
{\em $\texttt{[} 1, v \texttt{]}$} of the median group $(\wh{H}, \wh{m})$,
a bounded distributive lattice with respect to the partial order $\subset$,
is a deformation, induced by the median operation $m$ on $X$, of the
combinatorial configuration associated to the element $v$, in the following
sense. Let $C_v = \{w_i\,|\,i = \overline{1,\,n}\,\}, O_v = \{\zeta(w_i)\,|\,i = 
\overline{1,\,n}\,\}$, and the adjacent closed intervals $J_i, i = \overline{1,\,n}, n \geq 1$, 
as defined above. Then the following assertions hold.
\ben
\item[\rm (1)] $1 = \zeta(w_1) \subset \zeta(w_2) \subset \cdots
\subset \zeta(w_n) \subset v$.
\item[\rm (2)] The cell {\em $\texttt{[} 1, v \texttt{]}$} of the
median group $(\wh{H}, \wh{m})$ is the union of the adjacent
cells {\em $\wh{J}_i := \texttt{[} \zeta(w_i), \zeta(w_{i + 1})
\texttt{]} = w_i \texttt{[} \varphi(w_i^{- 1}), \varphi(w_i^{- 1} v)
\texttt{]} \sse w_i X$} for $i < n$ and {\em $\wh{J}_n :=
\texttt{[} \zeta(w_n), v \texttt{]} = w_n \texttt{[} \varphi(w_n^{- 1}),
\varphi(w_n^{- 1} v) \texttt{]} \sse w_n X$}, with the partial
order $\subset$ given by $u \subset u'$ for $u \in \wh{J}_i,
u' \in \wh{J}_k, i < k$, and {\em $u \subset u' \Llra w_i^{- 1} u \in
\texttt{[} \varphi(w_i^{- 1}), w_i^{- 1} u' \texttt{]} \sse X$}
for $u, u' \in \wh{J}_i$.
\een

Consequently, the median group $(\wh{H}, \wh{m})$ is
locally linear (simplicial) provided the median set $(X, m)$ is so.
Thus the restriction maps ${\rm res}_l : \cM_l(\wh{H}, \varphi) \lra
\cM_l(X)$ and ${\rm res}_s : \cM_s(\wh{H}, \varphi) \lra \cM_s(X)$
are bijective.
\end{co}

Thanks to Theorem 3.1. and Corollary 3.3., we can 
complete Lemma 1.14. as follows.

\begin{co}
Let $G$ be a group. Then the following assertions are equivalent.
\ben
\item[\rm (1)] $G$ acts freely on some nonempty median set (locally
linear median set, simplicial median set).
\item[\rm (2)] $G$ is embeddable into the underlying group of
some median group (locally linear median group, simplicial median
group).
\een
\end{co}

\begin{rem} \em
According to Lemma 1.14., the groups acting freely on median sets
form a quasivariety {\bf MSFG} axiomatized by very simple quasi-identities.
On the other hand, the class of groups acting freely on locally linear
median sets is axiomatized by the set of all universal sentences
in the first-order language of groups which are true in every
locally linear median group. It would be of some interest to find
a more concrete axiomatization, as well as a characterization of
the finitely generated members of the class above. Concerning
free actions on simplicial median sets, similar questions arise :
characterize the (finitely generated) groups acting freely on simplicial
median sets, as well as the models of the theory in the first-order
language of groups consisting of all universal sentences which are true 
in every simplicial median group.
\end{rem}

In particular, taking $H = 1$ in the statements above, we obtain

\begin{co}
Let $\X = (X, m)$ be a nonempty median set. Then there exist 
a median group $\F = (F, \wh{m})$ and an embedding $\iota : \X \lra \F$ of
median sets such that $1 \in \iota(X)$, $\iota(X) \sm \{1\}$ freely
generates the group $F$, and $\iota(X)$ is a retractible convex
subset of $\F$. In addition, the median group $\F$ is locally linear
(simplicial) provided the median set $\X$ is so.
\end{co}

\begin{ex} \em
Let $F$ be the free group of rank $3$ with generators $x_i, i = 1, 2, 3$,
and let $X = \{1, x_1, x_2, x_3\}$. Define the retract $\varphi : F \lra X$
of the embedding $X \lra F$ by 
\begin{center}
$\varphi(w) = \left\{ \begin{array}{lcl}
x_i  & \mbox{if} & x_i \leq w \\
1 & \mbox{if} & {\rm either}\, w = 1\, {\rm or}\, x_i^{- 1} \leq w\,{\rm for\, some}\,i \in \{1, 2, 3\}
\end{array}
\right.$
\end{center}
The set $\cM(X)$ of median operations on $X$ has $7$ elements : $4$
triangles and $3$ squares. If $m$ is the median operation of the triangle
with vertices $x_1, x_2, x_3$, so $m(x_1, x_2, x_3) = 1$, then the unique
median group operation $\wh{m} \in \cM(F, \varphi)$ with res$(\wh{m}) = m$
is the median operation of the canonical simplicial tree on $F$ determined
by the base $\{x_1, x_2, x_3\}$. If $m$ is the median operation of the
triangle with vertices $1, x_1, x_2$, so $m(1, x_1, x_2) = x_1 \cap x_2 = x_3$,
then its extension $\wh{m}$ is the median operation of the simplicial tree
on $F$ determined by the base $\{x_3^{- 1} x_1, x_3^{- 1} x_2, x_3\}$. The
median group operations corresponding to the other two triangles are obtained
in a similar way.

On the other hand, if $m$ is the median operation of the square with
the pairs $(1, x_2), (x_1, x_3)$ of opposite vertices, then the corresponding
median group $(F, \wh{m})$ is simplicial, but not locally linear,
 with the set $\{w \in F \sm \{1\}\,|\,\texttt{[} 1, w \texttt{]} = 
\{1, w\}\} = \{x_1^{\pm 1}, x_3^{\pm 1}, (x_1^{- 1} x_2)^{\pm 1}, (x_3^{- 1} x_2)^{\pm 1}\}$
of cardinality $8$, and the set $\{w \in F\,|\,\texttt{[} 1, w \texttt{]} = \square\} =
\{x_2^{\pm 1}, (x_1^{- 1} x_3)^{\pm 1}\}$ of cardinality $4$.
For any $v \in F \sm \{1\}$, the cell $\texttt{[} 1, v \texttt{]}$ is a finite
union of adjacent segments and/or squares. For instance, taking 
$v = x_1^{- 2} x_2 x_3^2 x_2$, the cell $\texttt{[} 1, v \texttt{]}$ is the
union of the $4$ adjacent segments $\texttt{[} 1, x_1^{- 1} \texttt{]} = x_1^{- 1}
{[} x_1, 1 \texttt{]}, \texttt{[} x_1^{- 1}, x_1^{- 1} x_2 \texttt{]} = 
x_1^{- 2} \texttt{[} x_1, x_2 \texttt{]}, \texttt{[} x_1^{- 2} x_2, x_1^{- 2} x_2 x_3
\texttt{]} = x_1^{- 2} x_2 \texttt{[} 1, x_3 \texttt{]}, \texttt{[} x_1^{- 2} x_2 x_3, 
x_1^{- 2} x_2 x_3^2 \texttt{]} = x_1^{- 2} x_2 x_3 \texttt{[} 1, x_3 \texttt{]}$ 
and of the square $\texttt{[} x_1^{- 2} x_2 x_3^2, v
\texttt{]} = x_1^{- 2} x_2 x_3^2 \texttt{[} 1, x_2 \texttt{]}$. The cases of the
other two squares on $X$, with the pairs of opposite vertices $(1, x_1), (x_2, x_3)$
and $(1, x_3), (x_1, x_2)$ respectively, are similar.

\end{ex}


\section{The relatively-transitive closure of a free action on a median set}

$\quad$ In this section we extend a given free action on a median set to
a larger one with a suitable universal property. By iterating this construction,
we shall obtain in the next section the {\em transitive closure} of
any given free action on a median set.

Let $H$ be a group acting freely on a median set $\X = (X, m)$. Fix as in
the previous sections a set $B = \{b_i\,|\,i \in I\}$ of representatives
of the $H$-orbits, with $1 \in I, I' := I \sm \{1\}$. Let $\wh{H} = H \ast F$
be the free product of $H$ and the free group $F$ with base $I'$, and identify
$X$ to the $H$-subset $H \sqcup (\sqcup_{i \in I'} H i) \sse \wh{H}$, with
the $H$-equivariant retract $\varphi : \wh{H} \lra X$.
With the notation above, the main result of this section (Theorem 2 from Introduction)
reads as follows.

\begin{te}
There exist a median set $\wh{\X} = (\wh{X}, \wh{\fm})$
and a free action of $\wh{H}$ on $\wh{\X}$ such that, identifying $\wh{H}$
with the $\wh{H}$-orbit of a base point of $\wh{\X}$, the composition of the maps
$X \lra \wh{H}, \wh{H} \lra \wh{X}$ is a $H$-equivariant embedding of median sets
$\X \lra \wh{\X}$ satisfying the following universal property.

$({\rm RTUP})$ For every group $\wt{H}$ acting freely on a median set
$\wt{\X} = (\wt{X}, \wt{m})$, every morphism $(\psi_0, \psi) : (H, \X) \lra
(\wt{H}, \wt{\X})$ in the category {\bf FAMS} of free actions on median sets,
satisfying $\psi(X) \sse \wt{H} \psi(1)$, extends uniquely to a morphism
$(\wh{\psi}_0, \wh{\psi}) : (\wh{H}, \wh{\X}) \lra (\wt{H}, \wt{\X})$ in
{\bf FAMS}.

In particular, the free action $(\wh{H}, \wh{\X})$ satisfying {\em (RTUP)}, 
called the {\em relatively-transitive closure} of the free action 
$(H, \X)$, is unique up to a unique isomorphism.
\end{te}

\bp
To construct the median set $\wh{\X}$, we consider the natural 
action from the right of $\wh{H}$ on the power set $\cP(\wh{H})$, 
$P^u := u^{- 1} P = \{u^{- 1} v\,|\,v \in P\}$ for 
$u \in \wh{H}, P \sse \wh{H}$, and let 
$$\fS := \{P \sse \wh{H}\,|\,\forall u \in \wh{H},\,P^u \cap X \in Spec\, \X\}.$$
The set $\fS$ is stable under the action of $\wh{H}$, is closed 
under the involution $P \mapsto \wh{H} \sm P$, and contains 
the $\wh{H}$-set $\cS := \{\varphi^{- 1}(\fp)^u\,|\,
\fp \in Spec\, \X, u \in \wh{H}\}$. Recall that, according to Theorem 3.1., 
$\cS = Spec\,(\wh{H}, \wh{m})$, where $\wh{m} : \wh{H}^3 \lra \wh{H}$ is 
the unique median group operation for which $\X = (X, m)$ is a retractible
convex median subset of $(\wh{H}, \wh{m})$ with associated folding $\varphi$.
Thus we obtain the $H$-equivariant embedding $Spec\,\X \lra \fS, \fp \mapsto
\varphi^{- 1}(\fp)$ with the $H$-equivariant retract $\fS \lra Spec\,\X,
P \mapsto P \cap X$.

$\wh{H}$ acts canonically on the power set $\cP(\fS)$ : $u\, \fM := \{P^{u^{- 1}} =
u P\,|\,P \in \fM\}$ for $u \in \wh{H}, \fM \sse \fS$, and the action is
compatible with the negation operator $\fM \mapsto \neg \fM := \{P \in \fS\,|\,
\wh{H} \sm P \notin \fM\}$. Consider the map $\wh{H} \lra \cP(\fS),
u \mapsto \fU(u) := \{P \in \fS\,|\,u \notin P\}$ which maps $\wh{H}$
onto the $\wh{H}$-orbit of $\fU(1)$. Composing the map above with the
restriction map $\cP(\fS) \lra \cP(\cS), \fM \mapsto \fM \cap \cS$,
we obtain the injective map $\wh{H} \lra \cP(\cS)$ (see the proof of Theorem 3.1.)
Consequently, the map $\wh{H} \lra \cP(\fS)$ is injective too, whence the
transitive action of $\wh{H}$ on the set $\{\fU(u)\,|\,u \in \wh{H}\}$ is free.
Notice also that $\neg \fU(u) = \fU(u)$ for all $u \in \wh{H}$.

Let us now consider the sublattice $\fL$ of the boolean algebra $\cP(\fS)$
generated by the subsets $\fU(u) = u\, \fU(1), u \in \wh{H}$. As $\es \in
\fU(u), \wh{H} \notin \fU(u)$ for all $u \in \wh{H}$, it follows that
$\fL \sse \cP(\fS) \sm \{\es, \fS\}$. More precisely, any element of $\fL$
has the form $\dbigcup_{i = 1}^n \fU(F_i),\, n \geq 1$, where
the $F_i$'s are nonempty  finite subsets of $\wh{H}$, and
$$\fU(F_i) := \bigcap_{u \in F_i} \fU(u) = \{P \in \fS\,|\,P \cap F_i = \es\},
\,i = \overline{1,\,n}.$$

Consequently, $\fL$ is closed under the negation operator $\neg$, so
$\fL$ is an unbounded distributive lattice with negation on which
$\wh{H}$ acts canonically. Moreover the subset 
$\wh{X} := \{\fM \in \fL\,|\,\neg\,\fM = \fM\}$ of the elements
of $\fL$, fixed by the negation operator $\neg$, is 
closed under the underlying median operation on the distributive lattice $\fL$
$$\fm(\fM_1, \fM_2, \fM_3) = (\fM_1 \cap \fM_2) \cup (\fM_2 \cap \fM_3) \cup
(\fM_3 \cap \fM_1) = (\fM_1 \cup \fM_2) \cap (\fM_2 \cup \fM_3) \cap (\fM_3 \cup \fM_1),$$
and the action of $\wh{H}$ on $\fL$ induces an action on the median set
$\wh{\X} = (\wh{X}, \wh{\fm})$.

The injective map $X \lra \wh{X}, x \mapsto \fU(x)$ identifies the
median set $\X = (X, m)$ with a median (not necessarily convex)
subset of $\wh{\X} = (\wh{X}, \fm)$, i.e.
$\fm(\fU(x), \fU(y), \fU(z)) = \fU(m(x, y, z))$ for $x, y, z \in X$.

According to Theorem 1.6., the restriction map 
$Spec\,\fL \lra Spec\,\wh{\X}, \fP \mapsto \fP \cap \wh{\X}$
maps homeomorphically the spectral space $Spec\,\fL$ of prime ideals
of the distributive lattice with negation $\fL$ onto
the spectral space $Spec\,\wh{\X}$, the dual of
the median set $\wh{\X}$. On the other hand, 
let us consider the $\wh{H}$-equivariant restriction map 
$$Spec\,\fL \lra \cP(\wh{H}), \fP \mapsto \fP \cap \wh{H} := 
\{u \in \wh{H}\,|\,\fU(u) \in \fP\}.$$

First let us show that for all $\fP \in Spec\,\fL$, $P := \fP \cap \wh{H} \in \fS$,
i.e. the map above takes values in $\fS$. We have to show that for all $u \in \wh{H}$,
$\fq := P^u \cap X \in Spec\,\X$. Let $x, y \in \fq$, i.e. $x, y \in X, u x, u y \in P$,
so $\fU(x), \fU(y) \in \fP^u := \{u^{- 1} \fM\,|\,\fM \in \fP\} \in Spec\,\fL$.
Since $\fU(z) \sse \fU(x) \cup \fU(y)$ for every element $z$ 
contained in the cell $\texttt{[} x, y \texttt{]}_\X$ of the median set $\X$,
and $\fP^u$ is an ideal of $\fL$, we deduce that $\fU(z) \in \fP^u$ for
all $z \in \texttt{[} x, y \texttt{]}_\X$, whence $\texttt{[} x, y \texttt{]}_\X \sse \fq$ 
for all $x, y \in \fq$. On the other hand, let $x, y \in X \sm \fq$, whence
$\fU(x), \fU(y) \notin \fP^u$. Assuming that $z \in \texttt{[}
x, y \texttt{]}_\X \cap \fq$, we obtain $\fU(x) \cap \fU(y) \sse \fU(z) \in \fP^u$.
As $\fP^u$ is a prime ideal of $\fL$, it follows that either $\fU(x) \in \fP^u$
or $\fU(y) \in \fP^u$, i.e. a contradiction.

Next let us show that the well defined restriction map $Spec\,\fL \lra \fS$ is
bijective. For any $P \in \fS$, we denote by Id$(P)$ the ideal of $\fL$
generated by the subset $\{\fU(u)\,|\,u \in P\}$.  Let us show that
Id$(P)$ is prime. Let $\fM_1, \fM_2 \in \fL$ be such that $\fM_1 \cap \fM_2 \in\,$
Id$(P)$, i.e. $\fM_1 \cap \fM_2 \sse \dbigcup_{i = 1}^n \fU(u_i)$ for some
$u_i \in P, i = \overline{1,\,n}, n \geq 1$. Obviously, 
$P \notin \fM_1 \cap \fM_2$, say $P \notin \fM_1$. As an element of 
$\fL$, $\fM_1 = \dbigcup_{j = 1}^m \fU(F_j)$
for some nonempty finite subsets $F_j \sse \wh{H}, j = \overline{1,\,m}$, 
with $m \geq 1$. Consequently, there exist 
$v_j \in F_j \cap P, j = \overline{1,\,m}$, therefore 
$\fM_1 \sse \dbigcup_{j = 1}^m \fU(v_j) \in\,$ Id$(P)$, whence
$\fM_1 \in\,$ Id$(P)$, i.e. Id$(P)$ is prime as desired. To conclude
that the map $Spec\,\fL \lra \fS$ is bijective, it remains to show that
$P =\,$ Id$(P) \cap \wh{H}$ for $P \in \fS$, and Id$(\fP \cap \wh{H}) = \fP$
for $\fP \in Spec\,\fL$. The inclusions $\sse$ are obvious, while the
opposite inclusions follow by straightforward verifications.

The bijection above identifies $\fS$  with $Spec\,\wh{\X} \cong Spec\,\fL$,
the dual of the median set $\wh{\X}$, generated (as median set) by
its subset $\{\fU(u)\,|\,u \in \wh{H}\}$, the $\wh{H}$-orbit of $\fU(1)$.
The $\wh{H}$-equivariant embedding $\cS = Spec\,(\wh{H}, \wh{m}) \lra \fS
\cong Spec\,\wh{\X}$ induces a $\wh{H}$-equivariant surjective 
morphism of median sets 
$\pi : \wh{\X} = (\wh{X}, \wh{\fm}) \lra (\wh{H}, \wh{m})$, with
$\pi|_X = 1_X$. In particular, it follows that the action of $\wh{H}$
on the median set $\wh{\X}$ is free. Notice also that we obtain a
canonical $H$-equivariant surjective morphism of median sets
$p : \wh{\X} \lra \X$, a retract of the $H$-equivariant embbeding
$\X \lra \wh{\X}$, by composing $\pi : \wh{X} \lra \wh{H}$ with
the folding $\varphi : \wh{H} \lra X$.

Finally, it remains to show that the free action $(\wh{H}, \wh{\X})$
extending $(H, \X)$ satisfies the universal property (RTUP). Let
$(\psi_0, \psi) : (H, \X) \lra (\wt{H}, \wt{\X})$ be a morphism in
{\bf FAMS} such that $\psi(X) \sse \wt{H} \psi(1)$, so we can identify
$\wt{H}$ with the $\wt{H}$-orbit $\wt{H} \psi(1) \sse \wt{X}$ and
$\psi(1) \in \wt{X}$ with the neutral element of the group $\wt{H}$,
whence the map $\psi : X \lra \wt{X}$ factorizes through $\wt{H}$
and $\psi_0 = \psi|_H$. As $X \sse \wh{H}$ is stable under the 
left multiplication with elements of $H$ and $\wh{H} = H \ast F$
is the free product of $H \sse X$ and the free group $F$ with base
$I' \sse X$, the map $\psi$ extends uniquely to the group morphism 
$\wh{\psi}_0 : \wh{H} \lra \wt{H}$ defined by
$\wh{\psi}_0(h) = \psi(h)$ for $h \in H$, $\wh{\psi}_0(i) = \psi(i)$
for $i \in I'$. To extend the group morphism $\wh{\psi}_0 : \wh{H} \lra
\wt{H}$ to the desired morphism $(\wh{\psi}_0, \wh{\psi}) : (\wh{H}, \wh{\X})
\lra (\wt{H}, \wt{\X})$ in {\bf FAMS} it suffices, by duality cf. Theorem 1.6.,
to define the morphism $Spec(\wh{\psi}) : Spec\,\wt{\X} \lra Spec\,\wh{\X} \cong \fS$ 
over $Spec\,\X$ in the unique possible way :
$\wt{P} \in Spec\,\wt{\X} \mapsto \wh{\psi}_0^{- 1}(\wt{P} \cap \wt{H})$.
The latter map is well defined, i.e. 
$\wh{\psi}_0^{- 1}(\wt{P} \cap \wt{H})^u \cap X = \psi^{- 1}(\wt{P}^{\wh{\psi}_0(u)}) 
= \psi^{- 1}(\wh{\psi}_0(u)^{- 1} \wt{P}) \in Spec\,\X$ for 
$u \in \wh{H}, \wt{P} \in Spec\,\wt{\X}$, since for all $x \in X$,
$u \in \wh{H}, \wt{P} \in Spec\,\wt{\X}$, we obtain
$$x \in \wh{\psi}_0^{- 1}(\wt{P} \cap \wt{H})^u \Llra \wh{\psi}_0(u x) = 
\wh{\psi}_0(u) \psi(x) \in \wt{P} \Llra$$
$$\psi(x) \in \wt{P}^{\wh{\psi}_0(u)} \in 
Spec\,\wt{\X} \Llra x \in \psi^{- 1}(\wt{P}^{\wh{\psi}_0(u)}) \in Spec\,\X.$$
Notice also that the map $Spec(\wh{\psi}) : Spec\,\wt{\X} \lra Spec\,\wh{\X}$
above is coherent, as required, since for every finite subset
$F \sse \wh{H}$, the inverse image of the basic quasicompact open 
set $\fU(F) = \{P \in \fS\,|\,P \cap F = \es\}$ of the spectral space
$\fS \cong Spec\,\wh{\X}$ is the quasicompact open set
$\{\wt{P} \in Spec\,\wt{\X}\,|\,\wt{P} \cap \wh{\psi}_0(F) = \es\}$ of
$Spec\,\wt{\X}$. This finishes the proof.
\ep

One checks easily that the correspondence $(H, \X) \mapsto (\wh{H}, \wh{\X})$
provided by the statement above extends to an endofunctor RTC : {\bf FAMS} $\lra$
{\bf FAMS} (the {\em relatively-transitive closure}), together with
a natural monomorphism rtc : $1_{\bf FAMS} \lra$ RTC.


\section{The transitive closure of a free action on a median set}

$\quad$ We are now in position to prove Theorem 3 and its median group
theoretic version Theorem 3' (see the Introduction).

Starting from the endofunctor RTC : {\bf FAMS} $\lra$ {\bf FAMS}
(the {\em relatively-transitive closure}) and the natural
monomorphism rtc : $1_{\bf FAMS} \lra$ RTC as defined in Section 4,
we consider the direct system $\R\T\C := ({\rm RTC}_n)_{n \in \N}$ of 
endofunctors of {\bf FAMS}, defined inductively by 
${\rm RTC}_0 = 1_{\bf FAMS}, {\rm RTC}_{n + 1} =
{\rm RTC} \circ {\rm RTC}_n$, with the connecting natural {\em monomorphisms}
${\rm rtc}_n : {\rm RTC}_n \lra {\rm RTC}_{n + 1}$, defined by 
${\rm rtc}_n = {\rm RTC}_n \cdot {\rm rtc}$ for $n \in \N$.

\begin{te}
The category {\bf FTAMS} of free and transitive actions on median sets
is a reflective full subcategory of {\bf FAMS}, i.e. the embedding 
$\iota :$ {\bf FTAMS} $\lra$ {\bf FAMS} has a left adjoint
{\em TC :} {\bf FAMS} $\lra$ {\bf FTAMS} (the {\em transitive
closure}). More precisely, the following assertions hold.
\ben
\item[\rm (1)] The endofunctor $\iota \circ {\rm TC}$ of {\bf FAMS}
is the {\em direct limit} of the direct system $\R\T\C$ of endofunctors
of {\bf FAMS}.
\item[\rm (2)] The natural transformation ${\rm TC} \circ \iota \lra
1_{\bf FTAMS}$, the {\em counit} of the ajunction, is a natural
{\em isomorphism}.
\item[\rm (3)] The natural transformation ${\rm tc}\, : 1_{\bf FAMS} \lra 
\iota \circ {\rm TC}$, the {\em unit} of the adjunction, is a natural
{\em monomorphism}.
\een
\end{te}

\bp
Let $(H, \X)$ be an object of {\bf FAMS}, so we may identify $H$
with the $H$-orbit of a base point of the median set $\X$ and the latter
with the neutral element $1$ of $H$. Applying step by step the 
endofunctor RTC : {\bf FAMS} $\lra$ {\bf FAMS},
we obtain a chain of suitable embeddings 
$$H_0 \lra X_0 \lra \cdots \lra H_n \lra X_n \lra H_{n + 1} \lra X_{n + 1}
\lra \cdots,$$ 
with $(H_0, \X_0) = (H, \X), (H_{n + 1}, \X_{n + 1}) = {\rm RTC}(H_n, \X_n)$ 
for $n \in \N$. It follows easily that the union ${\rm TC}(H, \X) :=
\dbigcup_{n \in \N} H_n = \dbigcup_{n \in \N} X_n$ becomes a median group,
and hence an object of {\bf FTAMS} extending $(H, \X)$ with the desired 
universal property (TUP) cf. Theorem 3  from Introduction. The rest of 
assertions concerning the adjunction between the categories {\bf FAMS} and
{\bf FTAMS} follow by straightforward verifications.
\ep




The next example ilustrates the complexity of the functorial construction 
above in the simplest nontrivial case of the free median group with one
generator.

\begin{ex} \em 
Let $\G = (G, m)$ be a nontrivial median group, and let $g \in G \sm \{1\}$.
The median subgroup $\wt{\G}$ of $\G$ generated by the element $g$ is
the union of the ascending chain $(G_n)_{n \in \N}$ of subgroups of $G$,
as well as the union of the ascending chain $(\X_n)_{n \in \N}$ of median
subsets of $(G, m)$, defined inductively by $G_0 = 1, X_0 = \{1, g\}$, 
$G_{n + 1}$ is the subgroup generated by $X_n \sse G$, while $\X_n$ is
the median subset generated by $G_n \sse G$ for $n \geq 1$. In other words, 
fhe free and transitive action of the at most countable group $\wt{G}$ 
on its underlying median set is the {\em transitive closure of} the trivial 
action of $G_0 = 1$ on $X_0 = \{1, g\}$ {\em inside} the free and 
transitive action of $G$ on its underlying median set. 
If $\G = \wt{\G}$ is the free
median group with one generator $g$, i.e. $\G \cong {\rm TC}(1, X_0)$, 
then $G_1 \cong \Z$, $G_n$ is a proper factor group of the free group
$G_{n + 1}$ for $n \geq 1$, so $G$ is a free group of countable rank, and
$\X_{n + 1}$ is the median set {\em freely generated by} $G_{n + 1} \sse X_{n + 1}$
{\em over} the median set $\X_n \sse G_{n + 1}$, for $n \in \N$, in particular, 
$\X_1 = {\rm fms}(\Z)$ is the free median set of countable rank.
\end{ex}

\end{document}